\documentclass[12pt,reqno]{amsart}
\usepackage[margin=1in]{geometry}
\usepackage{amsmath,amssymb,amsthm}
\usepackage{mathrsfs}
\usepackage{xcolor}
\usepackage{hyperref}
\allowdisplaybreaks

\numberwithin{equation}{section}
\newtheorem{Thm}{Theorem}[section]
\newtheorem{Def}{Definition}[section]
\newtheorem{Lem}{Lemma}[section]
\newtheorem{Pro}{Proposition}[section]
\newtheorem{example}{Example}[section]

\newtheorem{Ass}{Assumption}[section]
\newtheorem{Rem}{Remark}[section]

\renewcommand{\d}{\/\mathrm{d}\/}
\def\mR{\mathbb{R}}
 
\def\mN{\mathbb{N}} 
\def\W{\mathrm{W}} 
\def\F{\mathrm{F}}

\def\C{\mathrm{C}}
\def\D{\mathrm{D}}
\def\S{\mathrm{S}}
\def\mE{\mathbb{E}}

\def\mZ{\mathcal{Z}}
\def\Z{\mathrm{Z}}
\def\K{\mathrm{K}}

\def\Q{\mathrm{Q}}
\def\G{\mathrm{G}}

\def\mD{\mathcal{D}}
\def\mA{\mathscr{A}} 
\def\mcL{\mathscr{L}}

\def\mB{\mathscr{B}}
\def\mF{\mathscr{F}}
\def\mP{\mathbb{P}}
\def\mU{\mathscr{U}}
\def\A{{\bf A}} 
\def\B{{\mathrm B}}
\def\X{{\mathrm X}}
\def\Y{{\mathrm Y}}

\def\H{{ \mathrm{H}}} 
\def\L{ \mathrm{L}}
\def\V{ \mathrm{V}}
\def\A{ \mathrm{A}}

\def\U{{\mathrm U}}
\def\M{{\mathrm M}}
\def\ze{{\bf \zeta}}
\def\si{\sigma} 
\def\lam{\lambda}

\def\c{\wedge}

\def\ta{\tau_{l}}
\def\var{\varepsilon}
\def\b{\beta}
\def\al{\alpha}

\def\om{\omega}
\def\Om{\Omega}
\def\fr{\frac}

\def\wi{\widetilde } 
\def\tr{\mathrm{Tr}}

\def\la{\left(}
\def\ra{\right)}
\def\iZ{\int_{\mathcal Z}}

\def\i0t{\int_0^t}
\def\sk{\sum_{k=1}^\infty}

\let\originalleft\left
\let\originalright\right
\renewcommand{\left}{\mathopen{}\mathclose\bgroup\originalleft}
\renewcommand{\right}{\aftergroup\egroup\originalright}

\def\no{\nonumber}

\begin{document}
	
\newcommand{\Addresses}{{
		\bigskip
		\footnote{
	\noindent  \textsuperscript{1}Department of Mathematics, Indian Institute of Technology Roorkee-IIT Roorkee,
	Haridwar Highway, Roorkee, Uttarakhand 247 667, INDIA.
		\par\nopagebreak
	\noindent  \textit{e-mail:} \texttt{maniltmohan@gmail.com, manilfma@iitr.ac.in}
		
		\noindent \textsuperscript{2}Department of Mathematics,
		 Indian Institute of Space Science and Technology (IIST)
			Trivandrum- 695 547, INDIA. \par\nopagebreak \noindent
		\textit{e-mail:} \texttt{sakthivel@iist.ac.in}

		\textsuperscript{3} Go.AI, Inc, Beavercreek, OH 45431
		\par\nopagebreak \noindent
		\textit{e-mail:} \texttt{provostsritharan@gmail.com}
	
	$^*$Corresponding author.

}}}

\title[Dynamic Programming of  the Stochastic Burgers Equation]{Dynamic Programming of  the Stochastic Burgers Equation Driven by L\'evy Noise\Addresses}
\author[M. T. Mohan, K. Sakthivel  and S. S. Sritharan]
{ Manil T. Mohan\textsuperscript{1}, K. Sakthivel\textsuperscript{2}   and Sivaguru S. Sritharan\textsuperscript{3*}}
\maketitle
\begin{abstract}
	In this work, we study the optimal control of  stochastic Burgers equation perturbed  by Gaussian and L\'evy type noises with distributed  control process acting on the state equation. We use  the dynamic programming approach for the second order  Hamilton-Jacobi-Bellman (HJB) equation  consisting of an integro-differential operator with L\'evy measure  associated with the stochastic control problem. Using the regularizing properties of the transition semigroup corresponding to  the stochastic  Burgers equation and compactness arguments,  we solve the HJB equation and the resultant feedback control problem. 
\end{abstract}



\section{Introduction}
Optimal  control theory of fluid mechanics has been one of the important subjects in applied mathematics
with  several engineering  applications. An interesting problem in
this direction is the rigorous study of the feedback synthesis of optimal control problems
for the stochastic Burgers equation forced by  random noise  using the   infinite-dimensional Hamilton-Jacobi-Bellman
(HJB) equation associated with the problem. Noise term enters in the physical system as a forcing due to structural
vibration and other environmental effects that can be incorporated either as a random
boundary forcing or as a random distributed forcing in the state equation.  The optimal control of deterministic Navier-Stokes equations has been studied in \cite{Sr1} by showing that the value function, which is the minimum for an objective functional, is the  viscosity solution of the associated HJB equation and  the authors  \cite{G} extended this study  for  the optimal control of 2D stochastic   Navier-Stokes equations with Gaussian noise.    In \cite{D}, the authors considered the optimal control of stochastic Burgers equation with Gaussian noise. They introduced a specific control problem in which the covariance operator also acts on control process and using the Hopf transformation, they solved the  problem by  using dynamic programming approach. They removed this restriction in \cite{Da} and  solved  the control problem by again associating it with HJB equation. In the later case, they considered the mild form of the HJB equation and used the smoothing properties of the transition semigroup to obtain a smooth solution  to the HJB equation in a weighted space.   The same strategy has also been used to study the optimal control of 2D  stochastic Navier-Stokes equations with Gaussian noise in \cite{Da1}.   For various aspects of the optimal control of deterministic and stochastic fluid dynamic models, one may refer to \cite{Sr} and  also for other general systems, one can  look at the books \cite{Y,Ba}.

Unlike in the case of optimal control of stochastic Navier-Stokes/Burgers equations with Gaussian noise (see, \cite{Sr, Sri,G, D, Da, Da1}), there are only very few works in the case of control of stochastic PDEs with L\'evy noise (see, \cite{AA, MM}).  Moreover,  in the recent book \cite{BB},  a phenomenological study of fully developed turbulence  and intermittency is carried out and in which  it is nicely proposed that the experimental observations of these physical characteristics  can be modeled
by stochastic Navier-Stokes equations with L\'evy noise.  In fact, other  qualitative properties like, ergodicity and invariant measure for  stochastic Navier-Stokes/Burgers equations with L\'evy noise have already  been  studied in the literature (see, for example,\cite{Do, DX, MSS}).

In this paper, we  consider the optimal control of stochastic Burgers equation perturbed by Gaussian and L\'evy type noises with distributed stochastic control force acting on the equation.  Following the  methodology developed in \cite{Da}, we  transform the   HJB  equation of  partial integro-differential type with L\'evy measure into mild form using transition semigroup associated with the stochastic Burgers equation.  The regularity of  solutions for the  mild form of the HJB equation directly depends on the smoothing properties  of the semigroup and this is achieved via the Bismut-Elworthy-Li type formula {\cite{E}} derived for the Burgers equation with L\'evy noise. However,  boundedness of the derivative of the transition semigroup demands the finiteness of exponential moments of the stochastic Burgers equation. We have stated this exponential estimate which will be proved in the future work. Then the solution for the HJB equation is obtained in a space of smooth functions having weighted exponential growth by compactness arguments. This justifies the required smoothness of solutions for the feedback control formula and therefore by standard arguments, we prove the existence of an optimal pair for the control problem. 

This paper is organized as follows. In Section 2, we state the problem and give necessary function spaces to be used in the rest of the paper. In Section  3, we state the HJB equation and write the Galerkin approximation. The Section 4 is devoted to establish the crucial moment estimates for the Burgers equation with L\'evy noise. Using these estimates in Section 5, we prove  {\it a-priori} estimates  related to the smoothing properties of the semigroup. Finally in Section 6, we establish the main results of this paper. 

\section{Mathematical Formulation}
Let   $\U=\U(\xi,t)$  be a stochastic control process with values in the space of square integrable functions. Let $\Q$ be a positive, symmetric linear trace class operator. Let    $\W(\xi,t)$  and $\mathrm{N}(t,\cdot)$  be the Hilbert space valued cylindrical Wiener process  and  Poisson random measure respectively.       Consider the \emph{optimal control problem} for \emph{stochastic  Burgers equation} with distributed control    and Dirichlet boundary conditions:
\begin{equation}
\label{2p1}
\left\{
\begin{aligned}
\d\X(\xi,t)&=[\nu \D^2_\xi \X(\xi,t)+\D_\xi \X(\xi,t)\cdot \X(\xi,t)+\U(\xi,t)]\d t \\
&\quad+\Q^{1/2} \d\W(\xi,t) + \displaystyle \int_{\mZ}\G(\xi,t,z) \wi{\mathrm{N}}(\d t,\d z), \  (\xi,t)\in(0,1)\times (0,T) , \\
\X(0,t) &=\X(1,t) =0, t\in[0,T],\\
 \X(\xi,0)&=x(\xi), \xi\in[0,1],  
\end{aligned}
\right.
\end{equation}
where $\nu>0$ denotes the \emph{kinematic viscosity}.  We assume that  $\nu=1,$  and the processes $\W(t,\xi)$ and $\mathrm{N}(t,\cdot)$ are mutually independent. Though $\X,\U,\W$ and $\mathrm{N}$ are depending on the random parameter $\omega$,  for simplicity, we have not written  its dependence explicitly throughout the paper.

Since    the level of turbulence in a fluid flow can be characterized by the time averaged \emph{enstrophy},   it seems appropriate  to consider the minimization (over all controls from a suitable admissible set) of a  \emph{cost functional}
\begin{equation} \label{2p13a}
{\mathcal J}(0,T;x,\U)= \mE\left[\int_0^T\int_0^1\left(|\D_\xi \X(\xi,t)|^2+\frac{1}{2}|\U(\xi,t)|^2\right)\d\xi \d t+\int_0^1|\X(T,\xi)|^2\d\xi\right].
\end{equation}

\subsection{Functional  Setting}
We define the function spaces and notation frequently used in the sequel.  Let $\H=\L^2(0,1)$ be endowed with the inner product   $\left(\cdot,\cdot\right)$ and the norm $\|\cdot\|.$ Let $\V=\H_0^1(0,1)$ be endowed  with the norm $\left(\displaystyle\int_0^1 |\D_\xi\cdot|^2\d\xi\right)^{1/2}$  will be denoted as $\|\cdot\|_1$ to be consistent with the fractional powers. The induced duality, for instance between the spaces $\V$ and its dual $\V^\prime=\H^{-1}(0,1)$ will be denoted as $((\cdot,\cdot)).$   

Let us define the  positive self-adjoint operator $$\A u=-\D^2_\xi u\text{ for  }u\in \mD(\A)=\H^2(0,1)\cap \V.$$ Then  $\A^{-1}$ is a compact self-adjoint operator and hence by the spectral theory,  there exists sequence of  orthonormal basis functions $\{e_k\}_{k=1}^{\infty}$ in $\H$ and eigenvalues $\{\sigma_k\}$ accumulating at zero so that $\A^{-1}e_k=\sigma_k e_k, \ k\in\mN.$   Taking $\lambda_k=1/\sigma_k, $ we see that $\A e_k=\lambda_k e_k, \ k\in\mN,$  and 
$$0<\lam_1\leq \lam_2\leq \cdots\leq \lam_k\leq \cdots\to +\infty, \ \ \mbox{as} \ \ k\to +\infty.$$  
In this paper, we also use the fractional powers of $\A$.  For $u\in \H$ and  $\alpha>0,$ let us define
$$\A^\alpha u=\sum_{k=1}^\infty \lam_k^\alpha (u,e_k) e_k,  \ u\in\mD(\A^\alpha), $$ where  $$\mD(\A^\alpha)=\left\{u\in \H:\sum_{k=1}^\infty \lam_k^{2\alpha}|(u,e_k)|^2<+\infty\right\}.$$ 
Here  $\mD(\A^\alpha)$ is equipped with the norm 
\begin{equation} \label{fn}
\|\A^\alpha u\|=\left(\sum_{k=1}^\infty \lam_k^{2\alpha}|(u,e_k)|^2\right)^{1/2}.
 \end{equation}
Note that $\mD(\A^0)=\H,$  $\mD(\A^{\frac{1}{2}})=\V;$ in general, we set $\V_{\alpha}= \mD(\A^{\alpha/2})$ with $\|u\|_\alpha =\|\A^{\alpha/2} u\|.$   For any $s_1<s_2,$ the embedding $\mD(\A^{s_2})\subset \mD(\A^{s_1})$ is also compact. Indeed, it is known that  $\A$ has an orthonormal basis of eigenvectors $\{e_k\}$ given by $e_k(\xi)=\sqrt 2 \sin k\pi\xi, \ k\in\mN,  \ \xi\in[0,1]$ corresponding to the eigenvalues    $\lambda_k=k^2\pi^2, \ k\in \mN.$ Applying the H\"older inequality on the  expression \eqref{fn}, one can get the following interpolation estimate:
\begin{equation}\label{ie}
\|\A^{s}u\|\leq \|\A^{s_1}u\|^\theta\|\A^{s_2}u\|^{1-\theta},
\end{equation}
for any real $s_1\leq s\leq s_2$ and $\theta$ is given by $s=s_1\theta+s_2(1-\theta).$  

We define the nonlinear operator 
\[\B(u)=\frac{1}{2} \D_\xi(u^2), \; u\in \V.\]
Let us define the  trilinear operator 
\[b(u,v,w):\H\times \V\times \H\mapsto \mR, \ b(u,v,w)=\int_0^1 u(\xi) \D_{\xi}v(\xi) w(\xi) \d\xi. \]  An integration by parts yields 
\begin{equation}
b(u,u,v)=-\frac{1}{2}b(u,v,u)=b(v,u,u)  \ \ \mbox{and}  \ \  b(u,u,u)=0, \mbox{ for all } \; u,v\in \V.
\end{equation}

\subsection{L\'evy Noise and Hypothesis}
Let $(\Omega,\mF,\mP)$ be a complete probability space equipped with an increasing family of sub-sigma fields $\{\mF_t\}_{0\leq t\leq T}$ of $\mF$ satisfying usual conditions.

Since $\Q$ is a positive, symmetric and trace class operator on $\H,$  there exists  an orthonormal basis $\{f_k\}_{k=1}^{\infty}$ of $\H$  such that $\Q f_k=\varrho_kf_k, k\in \mathbb{N}.$ Here $\varrho_k$ is the eigenvalue corresponding to $\{f_k\}$ which is real and positive satisfying \[\tr(\Q)=\sk \varrho_k<+\infty \ \  \mbox{and} \ \ \Q^{1/2}v=\sk\sqrt{\varrho_k}(v,f_k)f_k, \ v\in \H.\]

The stochastic process $\{\W(t) : 0\leq t\leq T\}$ is an $\H$-valued  cylindrical Wiener process on $(\Omega,\mF,\{\mF_t\}_{t\geq 0},\mP)$  if and only if for arbitrary $t,$ the process $\W(t)$ can be expressed as  $\W(t)=\displaystyle\sum_{k=1}^\infty  \beta_k(t)f_k,$ where $\beta_k(t), k\in \mathbb{N}$ are independent, one dimensional Brownian motions on the space $(\Omega,\mF,\{\mF_t\}_{t\geq 0},\mP).$  


 Let $(\mZ,|\cdot|)$ be a separable Banach space and $({\bf L}_t)_{t\geq 0}$ be an $\mZ-$valued L\'evy process.  For every $\om\in\Om,$ ${\bf L}_t(\om)$ has at most countable number of jumps in an interval and the jump  $\Delta  {\bf L}_t(\om):[0,T]\to \mZ$ is  defined by $\Delta {\bf L}_t(\om):= {\bf L}_t(\om)-{\bf L}_{t-}(\om)$ at $t\geq 0.$  For $\mZ_0=\mZ\backslash\{0\},$ we define
\[\mathrm{N}([0,T],\Gamma)=\#\{t\in [0,T] : \Delta {\bf L}_t(\om)\in \Gamma\}, \  \mbox{where} \  \Gamma \in \mB(\mZ_0), \ \om\in\Om.\]
The measure $\mathrm{N}(\cdot,\cdot)$ is the Poisson random measure(or jump measure) with respect to $(\Omega,\mF,\{\mF_t\}_{0\leq t\leq T},\mP)$  associated with the L{\'e}vy process $({\bf L}_t)_{t\geq 0}.$ Here $\mB(\mZ_0)$ is the Borel  $\si$-field, $\mathrm{N}([0,T],\Gamma)$  is the Poisson random measure defined on $(([0,T]\times \mZ_0), \mB([0,T]\times\mZ_0)),$ and $\mu(\cdot)=\mE(\mathrm{N}(1,\cdot))$ is the $\sigma$-finite measure defined on $(\mZ_0, \mB(\mZ_0)).$ The intensity measure $\mu(\cdot)$ on $\mZ$ satisfies the conditions $\mu(\{0\})=0$   and 
\begin{equation}\label{2p10}
\int_{\mZ} \left(1\wedge |z|^p\right)  \mu(\d z)< +\infty,\ \ p\geq 2.
\end{equation}
The compensated Poisson random measure is defined by $\wi{\mathrm{N}}(\d t,\Gamma)=\mathrm{N}(\d t,\Gamma)- \mu(\Gamma)\d t,$ where $\mu(\Gamma)\d t$ is the compensator of the L\'evy process $({\bf L}_t)_{t\geq 0}$   and  $\d t$ is the Lebesgue  measure. 

Let  $\G:[0,1]\times[0,T]\times \mZ \to \H$ be a measurable and $\mF_t-$adapted process satisfying 
\[\mE\left[\left\|\int_0^T\int_{\mZ} \G(t,z)\wi{\mathrm{N}}(\d t,\d z)\right\|^2\right]<+\infty,\] where $\G(t,z):=\G(\xi,t,z)$. The  integral defined by  $\M(t):=\displaystyle\int_0^t\iZ \G(s,z)\wi{\mathrm{N}}(\d s,\d z)$ is an  $\H-$valued martingale and there exist an increasing c\`adl\`ag processes so-called \emph{quadratic variation process} $[\M]_t$ and \emph{Meyer process} $\langle \M\rangle_t$ such that $[\M]_t-\langle \M\rangle_t$ is a  martingale (see, \cite{Sa}).   
Moreover, we have the following It\^o isometry 
$$\mE\left[\left\|\int_0^T\iZ \G(t,z)\wi{\mathrm{N}}(\d t,\d z)\right\|^2\right]=\int_0^T\iZ   \|\G(t,z)\|^2\mu(\d z)\d t $$ and the relation $\mE\|\M(t)\|^2=\mE[\M]_t=\mE\langle \M\rangle_t.$ For more details on L\'evy processes one may refer to \cite{A} and \cite{P}.   Following are the basic assumptions on $\Q$, $\G(\cdot,\cdot)$ and $\mu(\cdot)$ used in the sequel. Other necessary assumptions on  the noise coefficient are stated  in relevant sections of the paper. 
\begin{Ass}\label{ass2.1}  
	\begin{itemize}  
	\item[ $(\A_1)$] For any $\kappa\in (1/2,1),$ the   operator $\Q$ satisfies  
	$$ \|\Q^{-1/2}x\|\leq C(\Q)\|\A^{\frac{\kappa}{2}}x\|  \  \  \mbox{for any} \  \  x\in \mD(\A^{\frac{\kappa}{2}}).$$ 
	\item[ $(\A_2)$] The jump noise coefficient  $\G(\cdot,\cdot)$ satisfies
	\begin{itemize}
	\item [(i)] $\displaystyle\sup_{t\in[0,T]}\int_{\mZ}\| \G(t,z)\|^{2} \mu(\d z)\leq C$.
		\item [(ii)] $\displaystyle\int_0^T\int_{\mZ}\| \G(t,z)\|_1^{p} \mu(\d z)\d t\leq C,$ for all $p\geq 2$.
	\item [(ii)$^{\prime}$]	The assumption (ii) and the Poincar\'e inequality give $\displaystyle\int_0^T\int_{\mZ}\| \G(t,z)\|^{p} \mu(\d z)\d t\leq C(\pi),$ for all $p\geq 2.$
	\end{itemize}
\end{itemize}
Also we fix measurable subsets $\Z_m$ of $\Z$ such that $\mu\left(\Z_m\right)<+\infty$ and
$\Z_m\uparrow \Z  \ \mbox{as} \ m\to\infty.$ 
\end{Ass}
In the rest of the paper, $C$ will denote a generic positive constant depending on the given arguments. 

\begin{example}
As an example of $\Q$, one may take $\Q=\A^{-\alpha}$, $\alpha\in(\frac{1}{2},1)$, since 
\begin{align*}
\tr(\Q)=\tr(\A^{-\alpha})=\sum_{k=1}^{\infty}\frac{1}{\lambda_k^{\alpha}}=\frac{1}{\pi^{2\alpha}}\sum_{k=1}^{\infty}\frac{1}{k^{2\alpha}}<+\infty,
\end{align*}
for $\alpha>\frac{1}{2}$. Since $\|\Q^{-1/2}x\|=\|\A^{\frac{\alpha}{2}}x\|,$ any $\alpha\in(\frac{1}{2},1)$ satisfies the Assumption $(\A_1)$.
	\end{example}

\subsection{Function Spaces}
To obtain the smoothing property of the transition semigroup, we need to work with some space of functions having exponential weights. 
Let us define  $B_R:=\{x\in \H: \|x\|\leq R\}$  and  for $k\in \mN,\al\in[0,1)$ the space
$$\C^{k+\al}(B_R):=\left\{ \psi:\H\to \mR\bigg| \psi\in \C^k(B_R) \ \mbox{and} \ \sup_{\substack{x,y\in B_R\\    x\neq y}} \frac{|\D^k\psi(x)-\D^k\psi(y)|_{\mathcal{L}^k(\H^k;\mR)}}{\|x-y\|^\al}<+\infty\right\},$$ where $\H^k=\underbrace{\H\times\cdots\times\H}_{k \text{ times }}$.
The space $\C^{k+\al}(B_R)$ is a Banach space with the norm  $$\|\psi\|_{\C^{k+\alpha}(B_R)}:=\sup_{x\in B_R} |\psi(x)|+\sum_{j=1}^k\sup_{x\in B_R}|\D^j\psi(x)|_{\mathcal{L}^j(\H^j;\mR)}+\sup_{\substack{x,y\in B_R\\    x\neq y}} \frac{|\D^k\psi(x)-\D^k\psi(y)|_{\mathcal{L}^k(\H^k;\mR)}}{\|x-y\|^\al}.$$
For $\varepsilon>0,$ we define the weighted space
$$\C_\var^{k+\al}(\H):=\left\{\psi:\H\to\mR\bigg|\ \psi\in \C^{k+\al}(B_R)  \text{ for all } R>0, \ \sup_{R>0} e^{-\var R^2}\|\psi\|_{\C^{k+\al}(B_R)}<+\infty\right\},$$ 
and 
\[\|\psi\|_{k+\al,\var}:= \sup_{R>0} e^{-\var R^2}\|\psi\|_{\C^{k+\al}(B_R)}, \;\psi\in \C_\var^{k+\al}(B_R).\]
For $\alpha=0$, we obtain the space $\C_{\var}^k$ as 
$$\C_\var^{k}(\H):=\left\{\psi:\H\to\mR\bigg|\ \psi\in \C^{k}(B_R)  \text{ for all } R>0, \ \sup_{R>0} e^{-\var R^2}\|\psi\|_{\C^{k}(B_R)}<+\infty\right\},$$ 
and 
\[\|\psi\|_{k,\var}:= \sup_{R>0} e^{-\var R^2}\|\psi\|_{\C^{k}(B_R)}, \;\psi\in \C_\var^{k}(B_R).\]
For $\var=0,$ we use  $\C^{k+\al}_0(B_R)$ and $\|\cdot\|_{k+\al,0}=\|\cdot\|_{k+\al}.$  Moreover,  we define \[\|\psi\|_{0,\var}=\sup_{R>0}e^{-\var R^2}\displaystyle\sup_{x\in B_R} |\psi(x)|,\]  and the following norms
$$\|\psi\|_{1,\var}=\|\psi\|_{0,\var}+\sup_{R>0}e^{-\var R^2}\sup_{h\in B_R}\frac{|\left(\D_x\psi(x),h\right)|}{\|h\|}$$
and 
$$\|\psi\|_{2,\var}=\|\psi\|_{0,\var}+\sup_{R>0}e^{-\var R^2}\displaystyle\sup_{x\in B_R} |\D_x\psi(x)|+\sup_{R>0}e^{-\var R^2}\sup_{h\in B_R}\frac{|\left(\D_
x^2\psi(x) h,h\right)|}{\|h\|^2}.$$

We wish to point out here that one can very well work with function spaces restricted to a ball of radius $R_{0}$ in the Hilbert space $\H$.

\section{The Hamilton-Jacobi-Bellman Equation}

The system is controlled through the process $\U:[0,1] \times [0,T]\times \Omega\to \H$ which is adapted to the filtration  $\{\mF_t\}_{t\geq 0}.$ For a fixed constant $\rho>0,$ we define the set of all admissible control  as 
\begin{equation}\label{3.0}
\mU^{0,T}_{\rho}=\Big\{\U\in \L^2(\Omega, \L^2(0,T;\H)):\|\U\|\leq\rho,\mP-\text{ a.s.,} \ \mbox{and}  \  \U \ \mbox{is adapted to} \{\mF_t\}_{t\geq 0}\Big\}.
\end{equation}
The  abstract form of the \emph{controlled stochastic  Burgers equation}  is given by
\begin{equation} \label{2p12}
\left\{\begin{aligned}
\d\X(t)&=[-\A\X(t)+\B(\X(t))+\U(t)]\d t+\Q^{1/2}\d\W(t)
+\displaystyle \int_{\mZ}\G(t,z) \wi{\mathrm{N}}(\d t,\d z),  \; t\in (0,T), \\
\X(0)&=x, \ x\in \H.
\end{aligned}
\right.
\end{equation}
The cost   functional  now associated with  the   optimal  control problem \eqref{2p12}  consists of the minimization over all controls $\U\in {\mathscr{U}_{\rho}^{0,T}} $ of the functional
\begin{equation} \label{2p13b}
{\mathcal J}(0,T;x,\U)= \mE\left[\int_0^T\left(\| \D_\xi \X(t)\|^2+\frac{1}{2}\|\U(t)\|^2\right)\d t+\|\X(T)\|^2\right].
\end{equation}
We aim to find a $\wi \U\in{\mathscr{U}_{\rho}^{0,T}} $ such that ${\mathcal J}(0,T;x,\wi \U)=\displaystyle{\inf_{\U\in {\mathscr{U}_{\rho}^{0,T}}} {\mathcal J}(0,T;x,\U).}$ 
First we  state the existence and uniqueness of solutions of the equation \eqref{2p12}. 
\begin{Def}
An $\H-$valued  c\`adl\`ag $\mF_t-$adapted process  $\X(t),t\in[0,T]$ is said to be a \emph{ solution} of \eqref{2p12} if 
$\B(\X(t))\in \H, $ a.s., and for arbitrary $t\in[0,T]$ and $\zeta\in \mD(\A),$
\begin{equation*}
\begin{aligned}
\left(\X(t),\zeta\right) =&\left( x,\zeta\right)-\int_0^t \left(\X(s),\A\zeta\right)\d s+\int_0^t\left(\B(\X(s))+\U(s),\zeta\right) \d s +\left( \Q^{1/2} \W(t),\zeta\right) \\
&+ \int_0^t \iZ \left( \G(t,z),\zeta\right) \wi{\mathrm{N}}(\d t,\d z), \ \ \mP\mbox{-a.s.}
\end{aligned}
\end{equation*}
In a given probability space $(\Omega,\mathscr{F},(\mathscr{F}_t)_{t\geq 0},\mathbb{P})$,  two solutions $\X_1(\cdot)$ and $\X_2(\cdot)$ are said to be \emph{pathwise unique} if $$\mathbb{P}\Big\{\omega\in\Omega:\X_1(t,\omega)=\X_2(t,\omega),\textrm{ for all }t\in[0,T]\Big\}=1.$$
\end{Def}

\begin{Thm}[Theorem 2.1, \cite{Do}] \label{thm3.1}Let $\G$ satisfy the assumption $(\A_2)(ii)'$  with $p=2,$  and  the control   $\U\in \L^2(\Om,\L^2(0,T;\H)).$ If $x\in \H,$ then there exists \emph{a unique  solution} $\X(\cdot)=\X(\cdot;x,\U)$ of the state equation \eqref{2p12} with trajectories in $\mathrm{L}^2(\Omega;\mathscr{D}([0,T];\H)\cap \L^2(0,T;\V)),$   where $\mathscr{D}(0,T;\H)$ is the space of all c\`adl\`ag functions from $[0,T]$ to $\H$. 
\end{Thm}

We use the \emph{dynamic programming approach} to  solve the optimal control problem which  involves the study of the \emph{value function}   defined as
${\mathscr V}(t,x):=\inf\limits_{\U\in {\mathscr{U}_{\rho}^{0,t}}} {\mathcal J}(0,t;x,\U),$  which is formally a  solution   of the following  \emph{Hamilton-Jacobi-Bellman equation} (see, Appendix \ref{DPP})
\begin{equation} \label{2p14}
\left\{
\begin{aligned}
\D_tv(t,x)&=\mcL_x v(t,x) +\phi(x)+ \displaystyle\inf_{\|\U\|\leq \rho}\left\{\left(  \U(t),\D_xv\right)+\frac{1}{2}\|\U(t)\|^2\right\},  \; t\in (0,T) , \\
  v(0,x)&= \|x\|^2, \ x\in \H,
  \end{aligned}
  \right.   
\end{equation}
where $\phi(x)=\| \D_\xi x\|^2$  and $\mcL_x v$ is \emph{the integro-differential  operator} given by
\begin{equation} \label{g1}
\begin{aligned}
\mcL_x v(t,x):=&\frac{1}{2}\tr(\Q \D_x^2v)+\left(-\A x+\B(x),\D_xv\right) \\
&+\iZ\big(v(t,x+\G(t,z))-v(t,x)-\left(\G(t,z),\D_xv\right)\big)\mu(\d z). 
\end{aligned}
\end{equation}
Note that   \emph{the Hamiltonian  function}  $$\F(p):= \inf_{\|\U\|\leq \rho}\left\{\left( \U(t),p\right)+\frac{1}{2}\|\U(t)\|^2\right\}$$ can be evaluated explicitly as 
\begin{equation} \label{2p15}
\F(p)=\left\{\begin{array}{lclclclc}
-\frac{1}{2}\|p\|^2   &\mbox{for}&  \| p\|\leq \rho, \\
-\rho\|p\|+\frac{\rho^2}{2}   &\mbox{for}&  \| p\|> \rho.
\end{array}\right.
\end{equation}
The HJB equation \eqref{2p14} can now be written as 
\begin{equation} \label{2p16}
\left\{
\begin{aligned}
\D_tv(t,x)=&\mcL_x v(t,x)+ \F(\D_x v(t,x))+ \phi(x), \ \ t\in (0,T), \\
 v(0,x)=& \|x\|^2, \ x\in \H.
  \end{aligned}
  \right.   
\end{equation}
Moreover, if $v$ is a smooth solution of the HJB equation \eqref{2p16}, the optimal control is  given by the formula  \begin{equation} \label{2p17}\wi \U(t)= {\mathcal G}( \D_xv(T-t,\wi \X(t))),\text{ where }
{\mathcal G}(p)=\left\{\begin{array}{lclclclc}
-p   &\mbox{for}&  \| p\|\leq \rho, \\
-\rho\frac{ p}{\| p\|}  &\mbox{for}&  \|p\|> \rho.
\end{array}\right.
\end{equation}
In (\ref{2p17}), $\wi\X(t)$ is the optimal solution of the following \emph{closed loop equation}
\begin{equation} \label{2p18}
\left\{\begin{aligned}
\d\wi \X(t)&=[-\A\wi \X(t)+\B(\wi \X(t))+{\mathcal G}(\D_xv(T-t,\wi \X(t)))]\d t+{\Q}^{1/2}\d\W(t)\\
&\quad+ \int_{\mZ}\G(t,z) \wi{\mathrm{N}}(\d t,\d z), \;t\in (0,T), \\
\wi  \X(0)&=x\in \H.
\end{aligned}
\right.
\end{equation}
The pair $(\wi  \X, \wi \U)$ is the \emph{optimal pair} of the control problem.   In order to obtain such a smooth solution to \eqref{2p16}, we use the transition semigroup  $(\S_t)_{t\geq 0}$ defined on    $(\Omega,\mathscr{F},\{\mF_s\}_{s\geq t},\mP)$    of the   uncontrolled stochastic Burgers equation  associated with \eqref{2p12}. The semigroup is defined by $(\S_tf)(x)=\mE[f(\Y(t,x))],$ where  $f\in \C_b(\H;\mR)$ and $\Y(t,x)$ is the solution of the following uncontrolled Burgers equation
\begin{equation}\label{2p19}
\left\{\begin{aligned}
\d\Y(t)&=[-\A\Y(t)+\B(\Y(t))]\d t+\Q^{1/2}\d\W(t)
\displaystyle+ \int_{\mZ}\G(t,z) \wi{\mathrm{N}}(\d t,\d z),   t\in (0,T), \\ \Y(0)&=x\in \H.
\end{aligned}
\right.
\end{equation}
Under the transition semigroup, the value function $v$ solves the following equation  in  mild form 
\begin{equation}  \label{2p19a}
v(t,x)=\S_tf(x)+\int_0^t\S_{t-s}\F(\D_xv(s,x))\d s+\int_0^t\S_{t-s}\phi(x)\d s,  \ \ \ \text{where} \ \ f(x)=\|x\|^2.
\end{equation}

\begin{Def}
	For any $x\in\H$, the function $v(\cdot,\cdot)$ is a \emph{mild solution} of (\ref{2p16}) if $v\in\C([0,T]\times\H)$, for any $t>0$, $v(t,\cdot)\in\C_{\tilde{\gamma}}^1(\H)$ for some $\tilde{\gamma}\in[0,\wi\var_0]$, where $\wi\var_0=\frac{\var_0}{2}$ defined in Proposition \ref{lem3.2}, and satisfies (\ref{2p19a}) for all $t\in[0,T]$. 
\end{Def}
As we are looking for a smoothness of the transition semigroup $(\S_t)_{t\geq 0}$, we appeal to \emph{Bismut-Elworthy-Li(BEL) formula} for the c\'adlag process which in turn requires the differential of $\Y(\cdot)$ corresponding to the system \eqref{2p19} with respect to the initial data. 

The first differential of $\Y(\cdot)$ in the direction of $h \in  \H$ is  defined by $\eta^h(t)=\left( \D_x\Y(t),h\right) $ and satisfies 
\begin{equation} \label{dd1}
\left\{\begin{aligned}
\D_t\eta^h(t)&=-\A\eta^h(t)+ \D_\xi\big(\Y(t)\eta^h(t)\big), \\
\eta^h(0)&=h.
\end{aligned}
\right.
\end{equation}
The second  differential of $\Y(\cdot)$ is    defined by $\zeta^{h}(t)= \left(\D^2_x\Y(t)h,h\right) $ and satisfies 
\begin{equation} \label{dd2}
\left\{\begin{aligned}
\D_t\zeta^{h}(t)&=-\A\zeta^{h}(t)+ \D_\xi\big(\eta^h(t)^2\big)+  \D_\xi\big(\Y(t)\zeta^h(t)\big),\\
\ze^h(0)&=0.
\end{aligned}
\right.
\end{equation}
All the equations and formulas stated above will be justified formally using approximation techniques.

\subsection{Finite Dimensional Approximation}  
We introduce the approximation for the controlled Burgers equation \eqref{2p12}.  Let $\{e_1,e_2,\cdots,e_m\}$ be the first $m$ eigenvectors of $\A$ and $\mathrm{P}_m$ be orthogonal projector of  $\H$ onto the space spanned by these $m$ eigenvectors. For a fixed $m\in \mN,$ define the approximation of the nonlinearity in the  state equation as
\begin{equation}
\label{3.13}
\B_m(x)=\frac{1}{2} \D_\xi g_m(x), \ \ x\in \mathrm{P}_m \V, \ \mbox{where} \ g_m(x)= \frac{mx^2}{m+x^2}, \ \ x\in \mR,\end{equation}
and the nonlinearity in the cost functional as (see section 3, \cite{Da})
\begin{equation}
\label{3.13a}
\phi_m(x)=\frac{m \|x\|_1^2}{m+\|x\|^2}, \ \  \ x\in \mathrm{P}_m\V \ \ \  \ \mbox{and} \ \ f_m(x)=\frac{m\|x\|^2}{m+\|x\|^2}, \ \ x\in \mathrm{P}_m\H.\end{equation}
\begin{Rem}\label{rem3.1}
	Let $f\in\mathrm{C}^1([0,1])$. For any $u\in\mathrm{H}_0^1(0,1)$, an integration by parts yields 
	\begin{align}\label{1}
	\left(\D_{\xi}f(u),u\right)_{\L^2}=\int_0^1\D_{\xi}f(u(\xi))u(\xi)\d \xi=-\int_0^1\underbrace{f(u(\xi))\D_{\xi}u(\xi)}_{\D_\xi(g(u(\xi)))}\d \xi,
	\end{align}
	where $g(u(\xi))=\displaystyle\int_0^{u(\xi)}f(r)\d r$, so that $g'(u(\xi))=f(u(\xi))$. From (\ref{1}), we also have
	\begin{align}\label{2}
	\left(\D_{\xi}f(u),u\right)_{\L^2}=-\int_0^1{\D_\xi(g(u(\xi)))}\d \xi=-g(u(\xi))\Big|_0^1=-g(u(1))+g(u(0))=0,
	\end{align} 
	since $u\in\H_0^1(0,1)$.
\end{Rem}
The Galerkin approximated  control system is 
\begin{equation} \label{2F1}
\left\{\begin{aligned}
\d  \X_m(t)&=[-\A  \X_m(t)+\mathrm{P}_m\B_m(\X_m(t))+\U_m(t)]\d t+\Q_m^{1/2}\d \W_m(t)\\
&\quad+\displaystyle \int_{\mZ_m}\G_m(t,z) \wi{\mathrm{N}}(\d t,\d z),  \ t\in (0,T), \\
\X_m(0)&=\mathrm{P}_m x, \ x\in \H,
\end{aligned}
\right.
\end{equation}
where $\U_m=\mathrm{P}_m \U, \ \Q_m=\mathrm{P}_m\Q $, $\W_m=\mathrm{P}_m\W$, $\G_m=\mathrm{P}_m\G$ and $\mathcal{Z}_m=\mathrm{P}_m\mZ$. In \eqref{2F1}, $\B_m(\cdot)$ is defined in \eqref{3.13}. The approximated optimal control problem is the minimization over  all $\U_m\in {\mathscr{U}_{\rho}^{0,T}} $ of the cost functional
\begin{equation} \label{2F2}
{\mathcal J_m}(0,T;\mathrm{P}_m x,\U_m)= \mE\left[\int_0^T\left(\phi_m(\X_m(t))+\frac{1}{2}\|\U_m(t)\|^2\right)\d t+f_m(\X_m(T))\right],
\end{equation}
where $\phi_m(\cdot)$ and $f_m(\cdot)$ are defined in \eqref{3.13a}. The finite dimensional  equation associated with \eqref{2p19} is given by 
\begin{equation} \label{2F3}
\left\{\begin{aligned}
\d\Y_m(t)&=[-\A\Y_m(t)+\mathrm{P}_m\B_m(\Y_m(t))]\d t+\Q_m^{1/2}\d\W_m(t) 
+ \int_{\mZ_m}\G_m(t,z) \wi{\mathrm{N}}(\d t,\d z), \   t\in (0,T),  \\
\Y_m(0)&=\mathrm{P}_m x, x\in \H.
\end{aligned}
\right.
\end{equation}
Then  by the semigroup definition $(\S^m_tf_m)(x)=\mE[f_m(\Y_m(t,x))],$  we have the following approximated  HJB equation:
\begin{equation} \label{2F4}
\left\{
\begin{aligned}
\D_tv_m(t,x)&=\mcL_x v_m(t,x)+ \F_m(\D_x v_m(t,x))+ \phi_m(x), \ \ t\in (0,T), \\
 v_m(0,x)&= f_m(x), \ x\in \mathrm{P}_m \H,
  \end{aligned}
  \right.
\end{equation}
where $\mcL_x v_m$ is the operator 
\begin{equation}\label{g2}
\begin{aligned}
\mcL_x v_m(t,x)&=\frac{1}{2}\tr(\Q_m \D_x^2v_m)+\left(-\A x+\B_m(x),\D_xv_m\right)  \\
&\quad+\int_{\mZ_m}\big(v_m(t,x+\G_m(t,z))-v_m(t,x)-\left(\G_m(t,z),\D_xv_m\right)\big)\mu(\d z). 
\end{aligned}
\end{equation}  
The mild form of the HJB equation \eqref{2F4} is
\begin{equation} \label{2p20}
v_m(t,x)=(\S^m_tf_m)(x)+\int_0^t\S^m_{t-s}\F_m(\D_xv_m(s,x))\d s+\int_0^t\S^m_{t-s}\phi_m(x)\d s.
\end{equation}
The approximation of \eqref{dd1} in the direction of $h\in \mathrm{P}_m\H$ is
\begin{equation}\label{d1}
\left\{
\begin{aligned}
\D_t\eta_m^h(t)
&=-\A\eta_m^h(t)+\frac{1}{2} \mathrm{P}_m\D_\xi\big(g_m^\prime(\Y_m(t))\eta_m^h(t)\big), \\
\eta_m^h(0)&=h
\end{aligned}
\right.
\end{equation}
and the approximation of the  second  differential of $\Y_m(\cdot)$  satisfies 
\begin{equation} \label{d2}
\left\{\begin{aligned}
\D_t\ze_m^{h}(t)&=-\A\ze_m^{h}(t)+\frac{1}{2} \mathrm{P}_m \D_\xi\big(g_m^{\prime\prime}(\Y_m(t))(\eta_m^h(t))^2\big)+ \frac{1}{2} \mathrm{P}_m\D_\xi\big(g_m^{\prime}(\Y_m(t))\zeta_m^h(t)\big), \\
\ze_m^h(0)&=0. 
\end{aligned}
\right.
\end{equation}
Now we state the following BEL formula   involving the   It\^o-L\'evy  process $\Y_m(\cdot)$ of the system \eqref{2F1}, whose proof is given in \cite{MSS}.
\begin{Pro}[Bismut-Elworthy-Li Formula]
For each $f\in \C_b(\mathrm{P}_m \H),$  the semigroup $(\S^m_t)_{t\geq0}$ is G\^ateaux differentiable  and  its derivative in any direction $h\in \mathrm{P}_m \H$ is given by 
\begin{equation} \label{BE}
\left(\D_x (\S^m_tf)(x),h\right)=\frac{1}{t}\mE\left[f(\Y_m(t,x))\int_0^t\left(\Q_m^{-1/2}\eta^h_m(s),\d\W_m(s)\right)\right], 
\end{equation}
for $x \in \mathrm{P}_m \H,  \ t\in[0,T].$  Moreover, the second differential is given by
\begin{equation} \label{BE6}
\begin{aligned}
\D^2_x (\S^m_tf)(x)\la h,h\ra&=\mE\left[f(\Y_m(t))\left\{\frac{1}{t}\int_0^t\la\Q_m^{-1/2}\zeta^h_m(s),\d\W_m(s)\ra\right\}\right]\\
&\quad+\mE\left[\la \D_x f(\Y_m(t)), \eta^h_m(t)\ra\frac{1}{t}\int_0^t\la\Q_m^{-1/2}\eta^h_m(s),\d\W_m(s)\ra\right].
\end{aligned}
\end{equation}
\end{Pro}
We use the following interpolation result frequently
\begin{Pro}[Lemma 4.2, \cite{Da}]\label{prop2.1}
	Let  $ \b_i\geq 0,$  $\gamma_i\geq 0,i=1,2$ with $\gamma=\lambda\gamma_1+ (1-\lambda)\gamma_2,\ \b=\lambda\b_1+ (1-\lambda)\b_2$ for any   $\lambda\in[0,1].$  If $\psi\in \C^{\b_1}_{\gamma_1}(\mathrm{P}_m\H)\cap \C^{\b_2}_{\gamma_2}(\mathrm{P}_m\H),$ then $\psi\in \C^{\b}_{\gamma}(\mathrm{P}_m\H)$ and there exists a constant $C>0$ depending on $\lambda,\b_i,\gamma_i,i=1,2$  satisfying  
	\begin{align} \label{I}
	\|\psi\|_{\b,\gamma} \leq C \|\psi\|^\lambda_{\b_1,\gamma_1} \|\psi\|^{ 1-\lambda}_{\b_2,\gamma_2}. 
	\end{align}
\end{Pro}

\section{ A-priori Estimates of the Stochastic Burgers Equation} 
In this section, we derive various moment estimates including energy estimates and exponential moment estimates for the stochastic Burgers equation with L\'evy noise. 
In particular, to estimate the first and second variation $\eta_m^h(\cdot)$ and $\zeta_m^h(\cdot)$ respectively (see Lemma \ref{lem3.3}) in \eqref{BE}  and \eqref{BE6}, we are in need of the following estimates.      

Let $\mathrm{K}(t):=\mathrm{I}(t)+\mathrm{J}(t)$, where 
$\mathrm{I}(\cdot)$ and $\mathrm{J}(\cdot)$ are the stochastic integrals defined by 
\begin{equation*}
\left\{
\begin{aligned}
\mathrm{I}(t)&:=\int_0^te^{-(t-s)\A}\Q^{1/2}\d\W(s),\\
\mathrm{J}(t)&:=\int_0^t\int_{{\mZ}}e^{-(t-s)\A}\G(s,z)\widetilde{\mathrm{N}}(\d s,\d z).
\end{aligned}
\right.
\end{equation*}
It is easy to show that $\mathrm{I}(\cdot)$ and $\mathrm{J}(\cdot)$ satisfy the following stochastic differential equations:
\begin{equation*}
\left\{
\begin{aligned}
\d\mathrm{I}(t)&=-\A\mathrm{I}(t)+\Q^{1/2}\d \W(t),\mathrm{I}(0)=0,\\ 
\d\mathrm{J}(t)&=-\A\mathrm{J}(t)+\int_{{\mZ}}\G(t,z)\widetilde{\mathrm{N}}(\d t,\d z),\mathrm{J}(0)=0.
\end{aligned}
\right.
\end{equation*}
\begin{Lem}[\cite{ADe,MSS}]\label{lem4.1a}
For any $T\geq 0$ and $p\geq 1$, we have the following:
\begin{align}
\mE\left[\sup_{t\in[0,T]}\|\mathrm{I}(t)\|_1^{2p}\right]&\leq C(p,T)[ \tr(\Q)]^p,\\  \mE\left[\sup_{t\in[0,T]}\|\mathrm{J}(t)\|_1^{2p}\right]&\leq C(p)\left[\left(\int_0^T\int_{\Z}\|\G(t,z)\|^2\lambda(\d z)\d t\right)^{p} +\int_0^T\int_{\Z}\|\G(t,z)\|^{2p}\lambda(\d z)\d t \right].
\end{align}
\end{Lem}

                                                                                                                                                                                                                                                                                                                                      \begin{Lem}\label{lem4.1}
Under the Assumption $(A_2)$, for any  $x\in\H$ and $p\geq 1$, the solution of (\ref{2F3}) satisfies the following a-priori estimate:
\begin{align}\label{ck}
\mE\left[\left(\sup_{t\in[0,T]}\|\Y_m(t)\|^2+\int_0^T\|\Y_m(t)\|_1^2\d t\right)^p\right]&\leq C(T,p,\Q)\left(1+\|x\|^{2p}\right).
\end{align} 
\end{Lem}   
\begin{proof}
Let us define $\mathrm{Z}_m(t)=\Y_m(t)-\mathrm{K}_m(t),$ where $\mathrm{K}_m(t)=\mathrm{P}_m\mathrm{I}(t)+\mathrm{P}_m\mathrm{J}(t)$. Then, $\mathrm{Z}_m(t)$ satisfies 
\begin{align}\label{3.31a}
\frac{\d\mathrm{Z}_m(t)}{\d t}=[-\A\mathrm{Z}_m(t)+\mathrm{P}_m\B_m(\mathrm{Z}_m(t)+\mathrm{K}_m(t))].
\end{align}
We multiply (\ref{3.31a}) with $\mathrm{Z}_m(\cdot)$ to get 
\begin{align}\label{3.32a}
\frac{1}{2}\frac{\d }{\d t}\|\mathrm{Z}_m(t)\|^2+\|\mathrm{Z}_m(t)\|_1^2&=(\B_m(\mathrm{Z}_m(t)+\mathrm{K}_m(t)),\mathrm{Z}_m(t)).
\end{align}
Using the Gagliardo-Nirenberg inequality, we have 
\begin{align}\label{sobo}
\|\Y_m\|_{\L^4(0,1)}^2\leq C\|\Y_m\|^{3/2}\|\Y_m\|_1^{1/2}.
\end{align} 
Using Remark \ref{rem3.1}, we have  $(\B_m(\Y_m),\Y_m)=0$, and H\"older's inequality,  (\ref{sobo}), and Young's inequality, we obtain 
\begin{align}\label{3.33a}
(\B_m(\mathrm{Z}_m+\mathrm{K}_m),\mathrm{Z}_m)&=(\B_m(\mathrm{Y}_m),\mathrm{Y}_m-\mathrm{K}_m)=-(\B_m(\mathrm{Y}_m),\mathrm{K}_m)\no\\&=-\frac{1}{2}\int_0^1\D_{\xi}\left(\frac{m\Y_m^2}{m+\Y_m^2}\right)\K_m\d \xi=-\int_0^1\frac{m^2\Y_m\D_{\xi}\Y_m}{(m+\Y_m^2)^2}\K_m\d\xi\no\\&\leq \|\Y_m\|_{\L^4(0,1)}\|\Y_m\|_1\|\K_m\|_{\L^4(0,1)}\leq C\|\Y_m\|^{3/4}\|\Y_m\|_1^{5/4}\|\K_m\|_{\L^4(0,1)}\no\\&\leq \frac{1}{4}\|\Y_m\|_1^2+C\|\K_m\|_{\L^4(0,1)}^{8/3}\|\Y_m\|^2\no\\&\leq \frac{1}{2}\|\Z_m\|_1^2+\frac{1}{2}\|\K_m\|_1^2+C\|\K_m\|_{\L^4(0,1)}^{8/3}\left(\|\Z_m\|^2+\|\K_m\|^2\right).
\end{align}
Substituting (\ref{3.33a}) in (\ref{3.32a}), we arrive at 
\begin{align}\label{3.34a}
\frac{\d }{\d t}\|\mathrm{Z}_m(t)\|^2+\|\mathrm{Z}_m(t)\|_1^2\leq \|\K_m(t)\|_1^2+C\|\K_m(t)\|_{\L^4(0,1)}^{8/3}\left(\|\Z_m(t)\|^2+\|\K_m(t)\|^2\right).
\end{align}
By the embeddings  $\H^1(0,1)\subset\subset\L^{\infty}(0,1)\hookrightarrow\L^4(0,1)\hookrightarrow\L^2(0,1)$, we also have
\begin{align}\label{3.34b}
\frac{\d }{\d t}\|\mathrm{Z}_m(t)\|^2+\|\mathrm{Z}_m(t)\|_1^2\leq \|\K_m(t)\|_1^2+C\|\K_m(t)\|_{1}^{8/3}\|\Z_m(t)\|^2+C\|\K_m(t)\|_{1}^{14/3}.
\end{align}
Integrating the above inequality from $0$ to $t$ to find 
\begin{align}
\|\mathrm{Z}_m(t)\|^2+\int_0^t\|\mathrm{Z}_m(s)\|_1^2\d s&\leq \|x\|^2+\int_0^t\|\K_m(s)\|_1^2\d s+C\int_0^t\|\K_m(s)\|_{1}^{14/3}\d s\no\\&\quad +C\int_0^t\|\K_m(s)\|_{1}^{8/3}\left(\|\Z_m(s)\|^2+\int_0^s\|\mathrm{Z}_m(r)\|_1^2\d r\right)\d s
\end{align}
An application of Gr\"onwall's inequality yields 
\begin{align}\label{3.40}
\|\Z_m(t)\|^2+\int_0^t\|\mathrm{Z}_m(s)\|_1^2\d s&\leq \left(\|x\|^2+\int_0^T\|\K_m(t)\|_1^2\d t+C\int_0^T\|\K_m(t)\|_{1}^{14/3}\d t\right)\no\\&\quad\times\exp\left(C\int_0^T\|\K_m(t)\|_{1}^{8/3}\d t\right)\no\\&\leq C(1+\|x\|^2)e^{CT}\exp\left(C\int_0^T\|\K_m(t)\|_{1}^{6}\d t\right),
\end{align}
for all $t\in[0,T]$. 
Taking expectation in \eqref{3.40} and applying Lemma \ref{lem4.1a}, we get
\begin{align}\label{3.41}
\mE\left[\left(\|\Z_m(t)\|^2+\int_0^t\|\mathrm{Z}_m(s)\|_1^2\d s\right)^p\right]&\leq C(p,T)(1+\|x\|^2)^p\exp\left(C(p,T)\mE\left[\sup_{t\in[0,T]}\|\K_m(t)\|_{1}^{6}\right]\right)\no\\&\leq C(p,T,\Q)(1+\|x\|^{2p}),
\end{align}
for all $t\in[0,T]$. Since $\|\Y_m(t)\|\leq \|\Z_m(t)\|+\|\K_m(t)\|$, one can complete the proof of the lemma. 
\end{proof}

\begin{Pro}\label{prop3.2}
Suppose that the Assumption $(\A_2)-(ii)^\prime$  is satisfied. Then for any  $p\geq 2,$ there exists a constant $C>0$  such that for any  $x\in \H,$ the following estimate holds:
\begin{equation} \label{2p31}
\begin{aligned}
\mE\left[\sup_{t\in[0,T]}\|\Y_m(t)\|^{p} \right]+& \mE\left[\int_0^T  \|\Y_m(t)\|^{p-2}\|\Y_m(t)\|_1^2 \d t \right] 
 \leq C(p, \Q,T)(1+\|x\|^{p}).  
\end{aligned}
\end{equation}
\end{Pro}

\begin{proof}
For a given $m\in \mN$ and for all $l>0,$ we define a sequence of stopping times 
\begin{align} \label{st1}
\tau_l:=\inf\Big\{t\leq T: \| \Y_m(t)\|\geq  l\Big\}. 
\end{align}
   The It\^o formula applied to $\|\Y_m(\cdot)\|^{p}$ leads to the following: 
\begin{equation} \label{2p32}
\mE\left[\sup_{r\in[0, t\c\ta]}\|\Y_m(r)\|^p\right]+p\mE\left[\int_0^{t\c\ta}\|\Y_m(s)\|^{p-2}\|\Y_m(s)\|^2_1 \d s\right]\leq\|x_m\|^p+\sum_{i=1}^5 J_i,
\end{equation}
where
\begin{align*}
J_1=&\mE\Bigg[\sup_{r\in[0, t\c\ta]}\bigg|\int_0^{r} \int_{\mZ_m}\Big(\|\Y_m(s-)+\G_m(s,z)\|^p-\|\Y_m(s-)\|^p\\
&-p\|\Y_m(s-)\|^{p-2}\la \G_m(s,z), \Y_m(s-)\ra\Big)\mathrm{N}(\d z,\d s)\bigg|\Bigg],\\
J_2=&p\mE\left[\sup_{r\in[0, t\c\ta]}\left|\int_0^r\|\Y_m(s)\|^{p-2}\la \Y_m(s), \Q_m^{1/2}\d\W_m(s)\ra\right|\right],  \\
J_3=&p\mE\left[\sup_{r\in[0, t\c\ta]}\left|\int_0^r \int_{\mZ_m} \|\Y_m(s-)\|^{p-2}\la \G_m(s,z),\Y_m(s-)\ra\wi{\mathrm{N}}(\d s,\d z)\right|\right],\\
J_4=&\frac{p(p-2)}{2}\mE\left[\int_0^{t\c\ta}\|\Y_m(s)\|^{p-4}\tr\big(\Q_m \Y_m(s)\otimes \Y_m(s)\big) \d s\right],\\
J_5=&\fr{p}{2}\mE\left[\int_0^{t\c\ta}\|\Y_m(s)\|^{p-2}\tr(\Q_m) \d s\right].
\end{align*}
Here we used the fact that  $\la \B_m(\Y_m(s)), \Y_m(s)\ra = 0$ (see Remark \ref{rem3.1}). 
For ${x},{y}\in \mathrm{P}_m\H$ and for any $p\geq 2,$   one  obtains from the Taylor formula that 
\begin{align} \label{2p28}
\big|\|{x}+{y}\|^p-\|{x}\|^p-p\|{x}\|^{p-2}({x},{y})\big|\leq C(p) (\|{x}\|^{p-2}\|{y}\|^2+\|{ y}\|^p).
\end{align}
Making use of the above inequality, we get 
\begin{align} \label{2p33}
J_1&\leq  \mE\bigg[\int_0^{t\c\ta} \int_{\mZ_m}\bigg(\bigg|\|\Y_m(s)+\G_m(s,z)\|^p-\|\Y_m(s)\|^p\no\\
&\quad-p\|\Y_m(s)\|_\al^{p-2}\la  \G_m(s,z), \Y_m(s)\ra\bigg|\bigg)\mu(\d z)\d s\bigg]\no\\
&\leq C(p)\mE\left[\int_0^{t\c\ta} \int_{\mZ_m}\Big(\|\Y_m(s)\|^{p-2}\|\G_m(s,z)\|^2+\|\G_m(s,z)\|^p\Big)\mu(\d z)\d s\right]\no \\ &\leq \frac{1}{8}\mE\left[\sup_{r\in[0,t\c\ta]}\|\Y_m(r)\|^{p}\right]+C(p)\mE\left[\left(\int_0^{t\wedge\ta}\int_{\mZ_m}\|\G_m(s,z)\|^2\mu(\d z)\d s\right)^{p/2}\right]\no\\&\quad+C(p)\mE\left[\int_0^{t\wedge\ta}\int_{\mZ_m}\|\G_m(s,z)\|^p\mu(\d z)\d s\right].
\end{align}
Using the Davis inequality  and Young's inequality, we estimate  $J_2$ as follows
\begin{align} \label{2p34}
J_2&\leq 3p \mE\left[\left(\int_0^{t\c\ta}\|\Y_m(s)\|^{2(p-1)} \tr(\Q_m) \d s\right)^{1/2}\right] \no \\ &\leq \mE\left[\sup_{r\in[0,t\c\ta]}\|\Y_m(r)\|^{p-1}\left(\int_0^{t\c\ta}\tr(\Q_m) \d s\right)^{1/2}\right]\no\\
&\leq\frac{1}{8} \mE\left[\sup_{r\in[0,t\c\ta]}\|\Y_m(r)\|^{p}\right]+C(p)T^{p/2}(\tr(\Q_m))^{p/2}.
\end{align}
 Applying the Burkholder-Davis-Gundy  type inequalities (see \cite{I}), we get
\begin{align} \label{2p35}
J_3&\leq3p\mE\left[\left(\int_0^{t\c\ta} \int_{\mZ_m} \|\Y_m(s)\|^{2(p-1)} \|\G_m(s,z)\|^2 \mu(\d z)\d s\right)^{1/2}\right] \no\\
&\leq 3p\mE\left[\sup_{r\in[0,t\c\ta]}\|\Y_m(r)\|^{p-1}\left(\int_0^{t\c\ta}\int_{\mZ_m} \|\G_m(s,z)\|^2 \mu(\d z)\d s\right)^{1/2}\right] \no \\
&\leq  \frac{1}{8} \mE\left[\sup_{r\in[0,t\c\ta]}\|\Y_m(r)\|^{p}\right]+C(p) \mE\left[\left(\int_0^{t\c\ta} \int_{\mZ_m} \|\G_m(s,z)\|^2 \mu(dz) \d s\right)^{p/2}\right],
\end{align}
where we also used Young's inequality in the last step.  Finally as above
\begin{align}\label{2p36}
J_4+J_5&\leq \frac{1}{8} \mE\left[\sup_{r\in[0,t\c\ta]}\|\Y_m(r)\|^{p}\right]+C(p)T^{p/2}(\tr(\Q_m))^{p/2}.
\end{align}
Moreover, it is easy to see that $ \tr(\Q_m)\leq \tr(\Q)$, $\|\G_m(\cdot,\cdot)\|\leq \|\G(\cdot,\cdot)\|$  and   $\|x_m\| \leq \|x\|.$
Substituting \eqref{2p33}-\eqref{2p36} back into \eqref{2p32}, we arrive at
\begin{align} \label{2p37}
&\mE\left[\sup_{r\in[0, t\c\ta]}\|\Y_m(r)\|^p\right]+2p\mE\left[\int_0^{t\c\ta}\|\Y_m(s)\|^{p-2}\|\Y_m(s)\|_1^2 \d s\right]\no\\
& \leq 2\|x\|^p+C(p)T^{p/2}(\tr(\Q))^{p/2}+ C(p)\left(\int_0^{T}\int_{{\mZ}}\|\G(s,z)\|^2\mu(\d z)\d\ s\right)^{p/2}\no\\&\quad+C(p)\int_0^{T}\int_{{\mZ}}\|\G(s,z)\|^p\mu(\d z)\d s.
\end{align}  
Let $$\Om_l := \Big\{\om\in \Om: \ \ \| \Y_m(t)\| <  l\Big\}.$$
From \eqref{2p37}, we get
\begin{align*} 
&\int_{\Om_l}\|\Y_m(t)\|^p \d\mP(\om) +\int_{\Om \backslash \Om_l}\|\Y_m(t)\|^p \d\mP(\om)
\leq C(p, \Q,T)(1+\|x\|^{p}).
\end{align*}
Since the first integrand is bounded, we get from the second integral that
$\mP (\Om \backslash \Om_l) \leq C/l^p.$  Moreover, $\mP(\om\in\Om: \    \ta <T) \leq \mP (\Om \backslash \Om_l) \leq C/l^p$ whence     $\lim\sup_{l\to \infty}\mP(\om\in\Om: \    \ta <T)=0.$   Thus,    $\ta\to T, \mP$-a.s. as $l\to\infty.$ By taking limit $l\to \infty$ in  \eqref{2p37} one can get the estimate \eqref{2p31} through Fatou's lemma . 
\end{proof}
\begin{Pro}\label{lem3.2} Let $\displaystyle 0<\var<\var_0$ be  small enough  and Assumption  \ref{ass2.1} hold. Then for any $x\in\H$ and  $t\in [0,T],$ the following holds:
		\begin{align}\label{e00}
			\mE\left[\exp\left(\var \|\Y_m(t)\|^2 + \var \int _0^t \|\Y_m(s)\|^2_1 \d s\right)\right] 
			\leq C(\var,\Q,\G,T)  \exp(\var \|x\|^2). 
		\end{align}
\end{Pro}
The proof follows easily for the case $x$ restricted to a ball of radius $R_{0}$ in $\H$ by an application of the Ito formula to exponential functions. This exponential estimate also holds up to a stopping time without restricting the initial data  $x$ to a ball in $\H$. Proof for these details are addressed in a future paper.

Next we derive the estimates concerning the differentials \eqref{d1} and \eqref{d2} using the exponential moment estimate stated in Proposition \ref{lem3.2}.

\begin{Lem}\label{lem3.3}
Suppose the conditions of Proposition \ref{lem3.2} hold true. For any $p\geq 2 $ and $x,h\in \mathrm{P}_m\H,$ there exists  a constant $C>0$  such that  for any $t\in [0,T]:$  
\begin{align}
\mE\left[\sup_{s\in [0,t]}\|\eta_m^h (s)\|^p +\i0t \|\eta_m^h(s)\|^{p-2}\|\eta_m^h(s)\|^2_1 \d s\right] &\leq C(\var,\Q,T,p) \ e^{\frac{\var p}{2}\|x\|^2} \|h\|^p,\label{4.32}\\
\mE\left[\left(\sup_{s\in [0,t]}\|\eta_m^h (s)\|^2 +\i0t \|\eta_m^h(s)\|^2_1 \d s\right)^p\right] &\leq C(\var,\Q,T,p) \ e^{\var p\|x\|^2} \|h\|^{2p},\label{4.33}
\end{align}
and
\begin{align} 
\mE\left[\sup_{s\in [0,t]}\|\zeta_m^h (s)\|^2 +\i0t \|\zeta_m^h(s)\|^2_1 \d s\right] &\leq C(\var,\Q,T) \ e^{\var\|x\|^2} \|h\|^4. 
\end{align}
Moreover, we have 
\begin{align}\label{3.77}
\sup_{\|h\|=1}\mE\left[\int_0^T\|\eta_m^{h}(s)-\eta^{h}(s)\|_1^2\d s\right]\to 0,\text{ as } m\to\infty.
\end{align}
\end{Lem}
\begin{proof} We multiply  \eqref{d1} by $p\|\eta_m^h(t)\|^{p-2}\eta_m^h(t)$ to obtain
\begin{align} \label{de1}
\lefteqn{\|\eta_m^h(t)\|^p +p \i0t \|\eta_m^h(s)\|^{p-2}\|\eta_m^h(s)\|_1^2 \d s }\no\\
&& = \|\mathrm{P}_m h\|^p + \frac{p}{2}\i0t \|\eta_m^h(s)\|^{p-2}  \la \mathrm{P}_m \D_\xi\big(g_m^\prime(\Y_m(s))\eta_m^h(s)\big), \eta_m^h(s)\ra \d s. 
\end{align}
Twice  integrating by parts with the note of $g_m^{\prime\prime}(x)\leq 2, \text{ for any }x\in \mR$ and using the Gagliordo-Nirenberg inequality,  we get 
\begin{align} \label{pe1}
&\|\eta_m^h\|^{p-2}\la \mathrm{P}_m \D_\xi\big(g_m^\prime(\Y_m)\eta_m^h\big), \eta_m^h\ra \no \\
&  = \|\eta_m^h\|^{p-2}\int_0^1 g_m^{\prime\prime}(\Y_m)(\D_\xi \Y_m) (\eta_m^h)^2 \d\xi \no\\
&\leq  \|\eta_m^h\|^{p-2}\|\Y_m\|_1 \|\eta_m^h\|^2_{\L^4(0,1)} \leq C \|\eta_m^h\|^{p-2}\|\Y_m\|_1 \|\eta_m^h\|^{1/2}_1 \|\eta_m^h\|^{3/2} \no\\
&\leq \frac{1}{2} \|\eta_m^h\|^{p-2}\|\eta_m^h\|_1^2 + C \|\Y_m\|_1^{4/3} \|\eta_m^h\|^p.
\end{align}
The  identity \eqref{de1} can be estimated by the Gr\"onwall inequality as
\begin{align}\label{de2}
\|\eta_m^h(t)\|^p + \frac{3p}{4}\i0t \|\eta_m^h(s)\|^{p-2}\|\eta_m^h(s)\|_1^2 \d s \leq  e^{C\i0t \|\Y_m(s)\|_1^{4/3} \d s }\|h\|^p.
\end{align}
For $\var>0,$ we get the following by Young's inequality:
\begin{align} \label{Y1}
 C\i0t \|\Y_m(s)\|_1^{4/3} \d s \leq C(\var,p,T)  + \frac{p\var}{2}\i0t \|\Y_m(s)\|_1^2 \d s,
 \end{align} 
whence
\begin{align}\label{3.81}
\mE\left[\sup_{s\in[0,t]}\|\eta_m^h(s)\|^p + \i0t \|\eta_m^h(s)\|^{p-2}\|\eta_m^h(s)\|_1^2 \d s\right] \leq e^{C(\var,p,T)} \mE\left(e^{\frac{\var p}{2}\i0t \|\Y_m(s)\|_1^2 \d s }\right)\|h\|^p.
\end{align}
Choosing $0<\var \leq \var_0,$ and applying the Proposition \ref{lem3.2}, one can complete \eqref{4.32}. Let us now take $p=2$ in \eqref{de2} to obtain 
\begin{align*}
\|\eta_m^h(t)\|^2 + \frac{3}{2}\i0t\|\eta_m^h(s)\|_1^2 \d s \leq  e^{C(\var,T)}e^{\var\i0t \|\Y_m(s)\|_1^{2} \d s }\|h\|^2.
\end{align*}
Thus it is immediate that 
\begin{align*}
\mE\left[\left(\sup_{s\in [0,t]}\|\eta_m^h(s)\|^2 + \i0t\|\eta_m^h(s)\|_1^2 \d s\right)^p\right] \leq  e^{pC(\var,T)}\mE\left[e^{p\var\i0t \|\Y_m(s)\|_1^{2} \d s }\right]\|h\|^{2p},
\end{align*}
and \eqref{4.33} follows easily from Proposition \ref{lem3.2}.

For the second part, we note that from \eqref{d2} that
\begin{align} \label{de3}
\|\zeta_m^h(t)\|^2 +2 \i0t \|\zeta_m^h(s)\|_1^2 ds 
 &=\i0t  \la \mathrm{P}_m \D_\xi\big(g_m^{\prime\prime}(\Y_m(s))(\eta_m^h(s))^2\big), \zeta_m^h(s)\ra \d s\no \\
 &\quad+  \i0t \la \mathrm{P}_m\D_\xi\big(g_m^{\prime}(\Y_m(s))\zeta_m^h(s)\big), \zeta_m^h(s)\ra \d s. 
\end{align} 
Integrating by parts and using $g_m^{\prime\prime}(x)\leq 2, x\in \mR$, we get
\begin{align*}
\left| \la \mathrm{P}_m \D_\xi\big(g_m^{\prime\prime}(\Y_m)(\eta_m^h)^2\big), \zeta_m^h\ra\right| 
\leq \frac{1}{2}\|\zeta_m^h\|_1^2+ 2 \|\eta_m^h\|^4_{\L^4(0,1)}.
\end{align*}
Note that by \eqref{de2} (with $p=6$ and $p=2$) and \eqref{Y1} (with weight $\frac{\var}{8}$), we get
\begin{align} \label{le4}
2\i0t \|\eta_m^h(s)\|^4_{\L^4(0,1)} \d s &\leq C\left(\i0t \|\eta_m^h(s)\|^6 \d s\right)^\frac{1}{2}  \left(\i0t \|\eta_m^h(s)\|^2_1 ds \right)^\frac{1}{2}\no \\
&\leq C \sqrt {t} e^{C(\var,T)} e^{\frac{\var}{2}\i0t \|\Y_m(s)\|_1^2 \d s }\|h\|^4 .
\end{align}
A calculation similar to \eqref{pe1} (with $p=2$) yields 
\begin{align*}\left|\la \mathrm{P}_m\D_\xi\big(g_m^{\prime}(\Y_m)\zeta_m^h\big), \zeta_m^h\ra\right| \leq C\|\Y_m\|_1^{4/3} \|\zeta_m^h\|^2 + \frac{1}{2}\|\zeta_m^h\|_1^2. \end{align*}
Taking \eqref{le4} into account  and applying the Gr\"onwall inequality, one can obtain
\begin{align} \label{de4}
\|\zeta_m^h(t)\|^2 + \i0t \|\zeta_m^h(s)\|_1^2 \d s 
 \leq C(T)e^{C(\var,T)} \Big(e^{C\i0t \|\Y_m(s)\|_1^{4/3} \d s}\Big)  \Big(e^{\frac{\var}{2}\i0t \|\Y_m(s)\|_1^2 \d s} \Big)\|h\|^4. 
\end{align}
Again using \eqref{Y1} and then Proposition \ref{lem3.2}, we can conclude the proof of second part.  The convergence (\ref{3.77}) can be established by similar arguments as in Lemma 3.3, \cite{Da}.
\end{proof}

Finally we prove  the almost sure convergence of the solutions of  \eqref{2F3}. 
\begin{Pro}\label{prop3.3}
	Suppose that the Assumption $(\A_2)(ii)$ is satisfied. Let $\{x_m\}_{m=1}^{\infty}$ be such that $x_m\to x$ strongly in $\H.$ Then $\Y_m(\cdot),$ the solution of  \eqref{2F3} converges to the unique solution  $\Y(\cdot)$ of \eqref{2p19}  in $\L^2(\Om;\mathrm{L}^{\infty}([0,T];\H))\cap   \L^2(\Om;\L^2(0,T;\V))$ having c\`adl\`ag paths
	and almost surely in $\mathscr{D}([0,T];\H)\cap\L^2(0,T;\V).$  
\end{Pro}
\begin{proof}
	The existence and uniqueness of solution of the  equation \eqref{2F3} is given in Theorem \ref{thm3.1}. We only need to prove the strong convergence of $\Y_m(\cdot)$ to $\Y(\cdot)$ in the given topology. Using the It\^o formula and taking expectation, we get
	\begin{align} \label{c1}
	\mE\left[\|\Y_m(t)\|^2+2\int_0^t\|\Y_m(s)\|^2_1\d s\right]
	&=\|x_m\|^2 +t\tr(\Q_m) +\int_0^t\int_{\mZ_m}\|\G_m(t,z)\|^2\mu(\d z) \d t\no\\&\leq \|x\|^2 +T\tr(\Q) +\int_0^T\iZ\|\G(t,z)\|^2\mu(\d z) \d t.
	\end{align}
	So, $\{\Y_m\}$ is bounded in $\L^2(\Om;\H)\cap   \L^2(\Om;\L^2(0,T;\V)).$  By the Banach-Alaoglu theorem, there exists a subsequence of $\{\Y_m\}$ which we denote again by $\{\Y_m\}$ such that 
	\begin{equation} \label{cc2}
	\left\{
	\begin{aligned}
	&\Y_m(t) \rightharpoonup \Y(t), \ \ \mbox{in} \ \ \L^2(\Om;\H), \\
	&\Y_m \rightharpoonup \Y, \ \ \mbox{in} \ \ \L^2(\Om;\L^2(0,T;\V)),  \ \mbox{as}  \ m\to\infty. 
	\end{aligned}
	\right.
	\end{equation}
	Applying the infinite dimensional It\^o formula (\cite{Me}) to $\|\Y(\cdot)\|^2,$ we get
	\begin{align} \label{c3}
	\mE\left[\|\Y(t)\|^2+2\int_0^t\|\Y(s)\|^2_1\d s\right]
	=\|x\|^2 +t\tr(\Q) +\int_0^t\iZ\|\G(s,z)\|^2\mu(\d z) \d s.
	\end{align}
	From (\ref{c1}) and (\ref{c3}), we obtain 
	\begin{align}
	&\mE\left[\|\Y_m(t)\|^2-\|\Y(t)\|^2+2\int_0^t\left(\|\Y_m(s)\|_1^2-\|\Y(s)\|_1^2\right)\d s\right]\no\\&\leq \|x_m-x\|^2+t\tr[(\mathrm{I}-\mathrm{P}_m)\Q)]+\int_0^t\int_{{\mZ}}\|(\mathrm{I}-\mathrm{P}_m)\mathrm{\G}(s,z)\|^2\mu(\d z)\d s\no\\&\leq \|x_m-x\|^2+t\|\mathrm{I}-\mathrm{P}_m\|_{\mathcal{L}(\H)}\tr(\Q)+\int_0^t\int_{{\mZ}}\|\mathrm{I}-\mathrm{P}_m\|_{\mathcal{L}(\H)}^2\|\mathrm{\G}(s,z)\|^2\mu(\d z)\d s\no\\& \to 0\text{ as }m\to\infty.
	\end{align}
	It leads to 
	\begin{align} \label{c4}
	\mE\left[\|\Y_m(t)\|^2+2\int_0^t\|\Y_m(s)\|^2_1\d s\right]
	\to \mE\left[\|\Y(t)\|^2+2\int_0^t\|\Y(s)\|^2_1\d s\right], \ \  \mbox{as} \ \ m\to \infty. 
	\end{align}
	Clearly, by the weak convergence of $\{\Y_m\}$ established in \eqref{cc2} together with \eqref{c4}, we get  the strong convergence of $\{\Y_m\}$ to $\Y$ in  $\L^2(\Om;\H)$ and $\L^2(\Om;\L^2(0,T;\V))$ respectively. Moreover, applying $\mathrm{P}_m$ on \eqref{2p19} and taking difference with \eqref{2F3}, we get
	\begin{equation}\d(\Y_m(t)-\mathrm{P}_m\Y(t))=-\A\big(\Y_m(t)-\mathrm{P}_m\Y(t)\big)\d t + \big(\mathrm{P}_m\B_m(\Y_m(t))-\mathrm{P}_m\B(\Y(t))\big)\d t.\end{equation}
	Taking inner product with $(\Y_m-\mathrm{P}_m\Y)$, integrating by parts and applying Young's inequality,  we obtain
	\begin{align}\label{3.68}
	&\frac{1}{2}\|\Y_m(t)-\mathrm{P}_m\Y(t)\|^2+\int_0^t\|\Y_m(s)-\mathrm{P}_m\Y(s)\|^2_1\d s\no\\
	&=\i0t \la  \mathrm{P}_m\B_m(\Y_m(s))-\mathrm{P}_m\B(\Y(s)),\Y_m(s)-\mathrm{P}_m\Y(s)\ra \d s\no\\& \leq \frac{1}{4}\int_0^t\|g_m(\Y_m(s))-\Y^2(s)\|^2\d s+\frac{1}{4}\int_0^t\|\Y_m(s)-\mathrm{P}_m\Y_m(s)\|_1^2\d s.
	\end{align}
	So, it leads to
	\begin{align}\label{3.69}
	\mE\left[\sup_{t\in[0,T]}\|\Y_m(t)-\mathrm{P}_m\Y(t)\|\right] &\leq \frac{1}{\sqrt{2}}\mE\left[\int_0^T \|g_m(\Y_m(s))-\Y^2(s)\|^2\d s\right]^\frac{1}{2}.
	\end{align}
	By definition of $g_m(x)$ given in  \eqref{3.13} ,we get
	\begin{align}\label{3.70}
	\|g_m(\Y_m)-\Y^2\|^2 &\leq 2\left\|\frac{m\Y_m^2}{m+\Y_m^2}-\Y_m^2\right\|^2+2\|\Y_m^2-\Y^2\|^2\no \\&\leq 2\left\|\frac{\Y_m^4}{m+\Y_m^2}\right\|^2+2\|\Y_m-\Y\|^2\|\Y_m+\Y\|^2.
	\end{align} 
	Using (\ref{3.70}) in (\ref{3.69}), we estimate the first integral 
	with the help of the inequality  $\displaystyle \frac{x^4}{m+x^2}\leq \frac{1}{2\sqrt{m}}x^3$, the Sobolev embedding $\H^{1/3}(0,1)\subset\L^6(0,1)$ and the interpolation inequality (\ref{ie}) as 
	\begin{align*}
	 \frac{1}{2\sqrt{m}}\mE\left[\int_0^T\|\Y_m(t)\|^6_{\L^6(0,1)}\d t\right]^{\frac{1}{2}}&\leq \frac{1}{2\sqrt{m}}\mE\left[\int_0^T\|\Y_m(t)\|^6_{\H^{1/3}(0,1)}\d t\right]^{\frac{1}{2}}\no\\&\leq \frac{1}{\sqrt{m}}\left(\mE\int_0^T\|\Y_m(t)\|^4\|\Y_m(t)\|_1^2\d t\right)^{\frac{1}{2}}\to0\text{ as }m\to\infty,
	\end{align*}
	by Proposition \ref{prop3.2}. Since  $\Y_m(\cdot)$ and $\mathrm{Y}(\cdot)$ are bounded in $\L^2(\Om;\L^\infty(0,T;\H))$ and in the first part, we proved that $\|\Y_m-\Y\|_{\L^2(\Om;\L^2(0,T;\V))}\to 0$ as $m\to\infty,$  the Poincar\'e inequality gives 
	\begin{align*}
	&\mE\left[\int_0^T   \|\Y_m(t)+\Y(t)\|^2\|\Y_m(t)-\Y(t)\|^2 \d t\right]^\frac{1}{2}\no \\
	&\leq C \left(\mE \left[\sup_{t\in[0,T]}(\|\Y_m(t)\|^2+\|\Y(t)\|^2)\right]\right)^\frac{1}{2} \left(\mE\left[\int_0^T\|\Y_m(t)-\Y(t)\|^2_1\d t\right]\right)^\frac{1}{2} \to 0\text{ as } m\to\infty.
	\end{align*}
	Therefore, for any fixed  $\varepsilon>0,$ we get by the Chebyshev inequality that
	\begin{align*}
	\mP\left(\sup_{t\in[0,T]}\|\Y_m(t)-\mathrm{P}_m\Y(t)\|\geq \varepsilon\right)&\leq \frac{1}{\varepsilon^2}\mE\left(\sup_{t\in[0,T]}\|\Y_m(t)-\mathrm{P}_m\Y(t)\|^2\right) \to 0 \  \mbox{as} \ \ m\to \infty
	\end{align*}
	and by the Markov inequality that 
	\begin{align*}
	\mP\left(\int_0^T\|\Y_m(t)-\mathrm{P}_m\Y(t)\|_1^2\d t\geq \varepsilon\right)\leq \frac{1}{\varepsilon}\mE\left(\int_0^T\|\Y_m(t)-\mathrm{P}_m\Y(t)\|_1^2\d t\right) \to 0 \  \mbox{as} \ \ m\to \infty.
	\end{align*}
	Then  there exists a subsequence of $\Y_m$ (denoted by $\Y_m$ for notational simplification, see Theorem 17.3, \cite{J}) such that 
	$\Y_m \to \Y$ almost surely in $\mathscr{D}([0,T];\H)\cap\L^2(0,T;\V).$ Since the solution $\Y(\cdot)$ is unique, the entire sequence $\Y_m(\cdot)$ converges almost surely. This completes  the proof.
\end{proof}

\section{Smoothing Properties of the Transition Semigroup} 
Next, we prove an estimate concerning the regularity of the  semigroup $\S_t^mf.$ The proof is similar to \cite{Da}.  In order to prove this estimate we are in need of the following estimate concerning  the stochastic integral in the BEL formula. For notational simplification, we suppress the dependence of $\Q,T,p$ in  constants appearing in the rest of the calculations. 
\begin{Lem}\label{lem3.4}
	Suppose the operator $\Q$ satisfies the assumption $(\A_1).$ Then for any $$\kappa\in \left(1/2,1\right),\; 2\leq p<\infty, \; h\in \H\text{ and }t\in [0,T],$$ there exists a constant $C>0$ such that 
	\begin{align}
	\Big(N_p(\eta_m^h)\Big)^{1/p} \leq C(\var) t^{\frac{(1-\kappa)}{2}} \ e^{\frac{\var}{2}  \|x\|^2} \|h\|. 
	\end{align} 
	where $\displaystyle N_p(\eta_m^h):=\mE\left[\left|\int_0^t\la\Q_m^{-1/2}\eta^h_m(s),\d\W(s)\ra\right|^p\right].$
\end{Lem}
\begin{proof}
	Applying the H\"older inequality and Burkholder-Davis-Gundy inequality (see Theorem 3.49, \cite{P}), we get
	\begin{align*}
	\Big(N_p(\eta_m^h)\Big) \leq  \mE\left(\int_0^t\|\Q_m^{-1/2}\eta^h_m(s)\|^2\d s\right)^{p/2}:=I.
	\end{align*}
	Since by $(\A_1),$ $\|\Q_m^{-1/2}\eta^h_m\|^2 \leq C(\Q) \|\A^{\fr{\kappa}{2}}\eta^h_m\|^2,$ we obtain again by the interpolation inequality \eqref{ie} (with $\theta=1-\kappa, s_1=0$ and $s_2=\frac{1}{2})$, H\"older's inequality and Lemma \ref{lem3.3} that
	\begin{align*}
	I &\leq  C\mE\left(\int_0^t\|\eta^h_m(s)\|^{2(1-\kappa)}\|\eta^h_m(s)\|_1^{2\kappa}\d s\right)^{p/2}\no\\
	&\leq C t^{\frac{(1-\kappa)p}{2}}\mE\left[\sup_{s\in[0,T]}\|\eta_m^h(s)\|^{(1-\kappa)p}\left(\int_0^t\|\eta_m(s)\|_1^2\d s\right)^{\frac{\kappa p}{2}}\right]\no\\
	&\leq C t^{\frac{(1-\kappa)p}{2}}\left(\mE\left[\sup_{s\in[0,t]}\|\eta^h_m(s)\|^{2(1-\kappa)p}\right]\right)^{\frac{1}{2}}\left[\mE\left(\int_0^t\|\eta^h_m(s)\|_1^2\d s\right)^{kp}\right]^{\frac{1}{2}}\no \\
	& \leq C(\var,\Q,T) t^{\frac{(1-\kappa)p}{2}} \ e^{\frac{\var p}{2} \|x\|^2} \|h\|^p. 
	\end{align*}
	This completes the proof. 
\end{proof}
\begin{Pro}\label{prop3.4}
Suppose Assumptions $(\A_1)$-$(\A_2)$ hold true. For any $\gamma<\wi\var_0=\frac{\var_0}{2},\al\in[0,1]$ and $\var>0,$ there exists a constant $C>0$ such that for any $f\in \C_\gamma^0(\mathrm{P}_m \H), \ t\in[0,T]:$ 
\begin{align} \label{s1}
\|\S_t^mf\|_{1+\al,\gamma+\var} \leq C(\alpha,\var,\gamma) t^{-(1+\al)(1+\kappa)/2}\|f\|_{0,\gamma}.
\end{align}
\end{Pro}
\begin{proof}  Let us begin with $\al=0.$ 
Since $\S_t^mf(x)=\mE\left[f(\Y_m(t))\big|\Y_m(0)=x\right],$ for any $\gamma<\var_0,$ we get by the conditional Jensen inequality, Markov property and Proposition \ref{lem3.2} that
\begin{align*}
|\S_t^mf(x)| &= \left|\mE\left[e^{\gamma\|\Y_m(t)\|^2}\left(e^{-\gamma\|\Y_m(t)\|^2}f(\Y_m(t))\right)\big|\Y_m(0)=x\right]\right|\no\\&\leq  \mE\left[e^{\gamma\|\Y_m(t)\|^2}\left(e^{-\gamma\|\Y_m(t)\|^2}\left|f(\Y_m(t))\right|\right)\big|\Y_m(0)=x\right]\no\\&\leq\sup_{R>0}\sup_{x\in B_R}e^{-\gamma\|x\|^2}|f(x)|\mE\left[e^{\gamma\|\Y_m(t)\|^2}\big|\Y_m(0)=x\right] \no\\ &\leq C(\gamma) \|f\|_{0,\gamma}e^{\gamma\|x\|^2},
\end{align*}
whence 
\begin{align} \label{s2}
\|\S_t^mf\|_{0,\gamma} \leq    C(\gamma)\|f\|_{0,\gamma}.
\end{align}
We first apply the H\"older inequality to  the BEL formula derived in \eqref{BE}. Then   for  $2\gamma  \leq \var_0,$    we obtain
\begin{align*}
\left|\la \D_x\S^m_tf(x),h\ra\right| &\leq 
\frac{1}{t}\mE\left(\big[e^{-\gamma\|\Y_m(t)\|^2}|f(\Y_m(t))|\big]e^{\gamma\|\Y_m(t)\|^2}\left|\int_0^t\la\Q_m^{-1/2}\eta^h_m(s),\d\W(s)\ra\right|\right)\no \\
&\leq \frac{1}{t} \|f\|_{0,\gamma} \left[\mE\left(e^{2\gamma  \|\Y_m(t)\|^2}\right)\right]^{1/2} \Big(N_2(\eta_m^h)\Big)^{1/2}\no\\
&\leq C(\gamma,\var)\frac{1}{t} \|f\|_{0,\gamma}  e^{\gamma\|x\|^2} \ t^{\frac{1-\kappa}{2}}e^{\frac{\var}{2} \|x\|^2} \|h\|,
\end{align*}
where $N_2(\eta_m^h)$ is defined in Lemma \ref{lem3.4} and the final estimate is due to Lemme \ref{lem3.2} and \ref{lem3.4}. Therefore, we get 
\begin{align} \label{s3}
\|\S_t^mf\|_{1,\gamma+\var} \leq    C(\gamma,\var) \ t^{-\frac{(1+\kappa)}{2}} \ \|f\|_{0,\gamma}.
\end{align}
Suppose $f\in \C_{\gamma+\var}^1(\mathrm{P}_m\H).$  Then the case $\al=1$ follows from  \eqref{BE6}.    Applying the H\"older inequality for two terms in \eqref{BE6} separately, we have  
\begin{align*}
|\la \D^2_x\S^m_tf(x) h,h\ra| &\leq  \frac{1}{t} \|f\|_{0,\gamma+\var} \left[\mE\left(e^{2(\gamma+\var)  \|\Y_m(t)\|^2}\right)\right]^{\frac{1}{2}} \Big(N_2(\zeta_m^h)\Big)^{\frac{1}{2}}\nonumber\\
&\quad+\frac{1}{t} \|f\|_{1,\gamma+\var} \left[\mE\left(e^{2(\gamma+\var) \|\Y_m(t)\|^2}\right)\right]^{1/2} \Big(\mE(\|\eta_m^h(t)\|^6)\Big)^{1/6}\Big(N_3(\eta_m^h)\Big)^{1/3}.
\end{align*}
For $0<2(\gamma+\var)\leq \var_0$, we have 
\begin{align*}
|\la \D^2_x\S^m_tf(x) h,h\ra| &\leq  C(\gamma,\var) \|f\|_{0,\gamma+\var} \  t^{-(1+\kappa)/2} \ e^{\left(\gamma+\frac{3\var}{2}\right)\|x\|^2}\|h\|^2 \no\\
&\quad +C(\gamma,\var) \|f\|_{1,\gamma+\var} \ t^{-(1+\kappa)/2}\ e^{(\gamma+2\var)\|x\|^2}\|h\|^2.
\end{align*}
It is easy to deduce the following
\begin{align*}
\|\S_t^mf\|_{2,\gamma+2\var} \leq    C(\gamma,\var)  \ t^{-(1+\kappa)/2} \ \|f\|_{1,\gamma+\var}.
\end{align*}
By the semigroup property and \eqref{s3}, we get
\begin{align} \label{s4}
\|\S_t^mf\|_{2,\gamma+2\var} &= \|\S_{t/2}^m(\S_{t/2}^mf)\|_{2,\gamma+2\var}  \leq   C(\gamma,\var)  t^{-\frac{(1+\kappa)}{2}} \ \|\S_{t/2}^mf\|_{1,\gamma+\var}\no\\
&\leq C(\gamma,\var)  t^{-(1+\kappa)} \ \|f\|_{0,\gamma}.
\end{align}
Since $0<\var<\var_0$ (see Proposition \ref{lem3.2}) is arbitrary and so the result follows for $\al=1.$ For $\al\in(0,1),$ one can use the interpolation estimate given in Proposition \ref{prop2.1}.   Taking  $\lambda=1-\alpha$ and choosing $s_1=1,s_2=2,$ we get  
\begin{align} \label{s5}
\|\S_t^mf\|_{1+\alpha,\gamma+\var} &\leq C(\gamma,\var)   \|\S_t^mf\|^{1-\alpha}_{1,\gamma+\var}\|\S_t^mf\|^{\alpha}_{2,\gamma+\var}\leq C(\gamma,\var)  \ t^{-\frac{(1+\kappa)(1+\alpha)}{2}}  \|f\|_{0,\gamma}.
\end{align}
This completes the proof.
\end{proof}
In order to estimate the mild form of the HJB equation \eqref{2p20}, we need the following estimate.

\begin{Lem}\label{lem3.5}
Let $\al\in[0,1].$ For any $\var>0, \gamma<\var_0,$ there exists a constant $C>0$ such that  for any $t\in[0,T]$ and $m\in\mN,$ we have
\begin{align} \label{i1}
\left\|\int_0^t\S^m_{t-s}\phi_m( x)\d s\right\|_{1+\al,\gamma+\var} \leq C(\al,\gamma,\var).
\end{align}
\end{Lem}
\begin{proof}
For $x,h\in \mathrm{P}_m\V,$ we obtain
\begin{align*}
\la \D_x\int_0^t\S^m_{t-s}\phi_m( x)\d s,h\ra =  \mE\left[\int_0^t\la \D_x\phi_m(\Y_m(t-s)),\eta_m^h(t-s)\ra \d s\right]=:I_1.
\end{align*}
Note that 
\begin{align}\label{st}
\left(\D_x\phi_m(x),\eta_m^h\right)&=\left(\D_x\left(\frac{m\|x\|_1^2}{m+\|x\|^2}\right),\eta_m^h\right)\leq |(\D_x\|x\|_1^2,\eta_m^h)|+\|x\|_1^2|(\D_x\|x\|^2,\eta_m^h)|\no\\&\leq 2\|x\|_1\|\eta_m^h\|_1+2\|x\|_1^2\|x\|\|\eta_m^h\|.
\end{align}
Thus by using H\"older's inequality, Lemma \ref{lem4.1}, Proposition \ref{prop3.2} and Lemma \ref{lem3.3}, we have 
\begin{align*}
I_1&\leq 2\left[\mE\left(\int_0^t\|\Y_m(s)\|_1^2\d s\right)\right]^{1/2}\left[\mE\left(\int_0^t\|\eta_m^h(s)\|^2_1\d s\right)\right]^{1/2}\no\\&\quad +2\left[\mE\left(\sup_{s\in[0,t]}\|\Y_m(s)\|^4\right)\right]^{1/4}\left[\mE\left(\sup_{s\in[0,t]}\|\eta_m^h(s)\|^4\right)\right]^{1/4}\left\{\mE\left[\left(\int_0^t\|\Y_m(s)\|_1^2\d s\right)^2\right]\right\}^{1/2}\no\\&\leq C(\var)(1+\|x\|^2)^{1/2}e^{\frac{\var}{2}\|x\|^2}\|h\|+C(\var)(1+\|x\|^4)^{3/4}e^{\frac{\var}{2}\|x\|^2}\|h\|\no\\&\leq  C(\gamma,\var)e^{(\gamma+\var)\|x\|^2}\|h\|.
\end{align*}
Hence we get the result for $\al=0$ as
\begin{align*}
\left\|\int_0^t\S^m_{t-s}\phi_m(x) \d s\right\|_{1,\gamma+\var} \leq C(\gamma,\var).
\end{align*}
For the case of $\al=1,$ we note that 
\begin{align*}
&\la \D^2_x\left[\int_0^t\S^m_{t-s}\phi_m(x)\d s\right] h,h\ra \\&=  \mE\left[\int_0^t \big(\la \D_x\phi_m(\Y_m(t-s)),\zeta_m^h(t-s)\ra
+\D^2_x\phi_m(\Y_m(t-s))\la\eta_m^h(t-s),\eta_m^h(t-s)\ra \big)\d s\right]\\& =:I_2+I_3.
\end{align*}
A standard computation similar to \eqref{st} leads to the following
\begin{align*}
I_3&\leq 8\mE\left[\int_0^t\|\Y_m(s)\|\|\eta_m^h(s)\|\|\Y_m(s)\|_1\|\eta_m^h(s)\|_1\d s\right] +8\mE\left[\int_0^t\|\Y_m(s)\|^2\|\eta_m^h(s)\|^2\|\Y_m(s)\|_1^2\d s\right]\no\\&\quad +2\mE\left[\int_0^t\|\eta_m^h(s)\|^2\|\Y_m(s)\|_1^2\d s\right]+2\mE\left[\int_0^t\|\eta_m^h(s)\|_1^2\d s\right]\no \\&\leq C(\gamma,\var)e^{(\gamma+\var)\|x\|^2}\|h\|^2.
\end{align*}
A computation similar to $I_1$ for $I_2$ and the estimate of $I_3$ give
\begin{align*}
\left\|\int_0^t\S^m_{t-s}\phi_m(x)\d s\right\|_{2,\gamma+\var} \leq C(\gamma,\var).
\end{align*}
The other cases of $\al\in(0,1)$ can be done, as we argued in Proposition \ref{prop3.4}, by using the interpolation result in Proposition \ref{prop2.1}. 
\end{proof}

Using the smoothing property of the transition semigroup proved in Proposition \ref{prop3.4} and Lemma \ref{lem3.5}, we prove the regularity of solutions for the mild form \eqref{2p20} of the HJB equation \eqref{2F4}.
\begin{Pro}\label{prop3.5}
Let $\al\in(0,1)$ be fixed such that $(1+\al)(1+\kappa)<2$ and let $\gamma<\wi\var_0.$  Suppose the Assumptions $(\A_1)-(\A_2)$ hold true. Then for any $\var>0,$ there exist a constant $C>0$ such that the following holds:
\begin{align}\label{Dv}
\sup_{t\in[0,T]}\Big(t^{(1+\al)(1+\kappa)/2}\|v_m(t,\cdot)\|_{1+\al,\gamma+\var}\Big)\leq C(\al,\gamma,\var).
\end{align}
and  
\begin{align}\label{DDv}
\sup_{t\in[0,T]}\Big(t^{(1+\kappa)}\|v_m(t,\cdot)\|_{2,\gamma+\var}\Big)\leq C(\gamma,\var).
\end{align}
\end{Pro}
\begin{proof}
Note that the solution $v_m(t,x)$ which solves the HJB equation written in mild form \eqref{2p20} is in fact the value function to the approximate control problem on $[0,T]$:
\begin{align*}
v_m(T,x)&=\inf_{\U_m\in\mathscr{U}_\rho^{0,T}} \mathcal{J}_m(0,T;\mathrm{P}_mx,\U_m) \leq \mathcal{J}_m(0,T;\mathrm{P}_mx,0) \\
&=  \mE\left[\int_0^T\left(\frac{m\|\Y_m(t)\|^2_1}{m+\|\Y_m(t)\|^2}\right)\d t+\frac{1}{2} \frac{m\|\Y_m(T)\|^2}{m+\|\Y_m(T)\|^2}\right]\\
&\leq  \mE\left[\int_0^T\|\Y_m(t)\|^2_1\d t\right]+ \mE\left[\|\Y_m(T)\|^2\right]\leq C(\Q,T)(1+\|x\|^2).
\end{align*}
The last estimate follows from Proposition \ref{prop3.2}.  Hence, we get 
\[\|v_m(T,\cdot)\|_{0,\gamma}\leq C(\gamma).\]
But, if we consider the control problem in $[0,t],$ for any $t\in [0,T],$ one can also get
\begin{align}\label{ve1}
\|v_m(t,\cdot)\|_{0,\gamma}\leq C(\gamma).
\end{align}
From \eqref{2p20}, Proposition \ref{prop3.4} and Lemma \ref{lem3.5}, we obtain
\begin{align}  \label{ve2}
&\|v_m(t,\cdot)\|_{1+\al,\gamma+\var}\no\\&\leq\|\S^m_tf_m\|_{1+\al,\gamma+\var}+\int_0^t\|\S^m_{t-s}\F_m(\D_xv_m(s))\|_{1+\al,\gamma+\var}\d s+\left\|\int_0^t\S^m_{t-s}\phi_m\d s\right\|_{1+\al,\gamma+\var} \no\\
&\leq
 C(\al,\gamma,\var)(1+t^{-(1+\al)(1+\kappa)/2}\|f\|_{0,\gamma}) +\int_0^t\|\S^m_{t-s}\F_m(\D_xv_m(s))\|_{1+\al,\gamma+\var}\d s.
\end{align}  
The final term in (\ref{ve2}) can be estimated using the interpolation result of Proposition \ref{prop2.1}, Proposition \ref{prop3.4} and \eqref{ve1} and arguing similarly as in Proposition 4.4 of \cite{Da}, for any $\var<(\wi\var_0-\gamma)(1+\alpha)$.   The estimate \eqref{DDv} follows from Proposition 4.5 of \cite{Da}. 
\end{proof}

In order to prove the convergence of $v_m(\cdot,\cdot)$ in the mild form (\ref{2p20}), we need more spatial regularity on $v_m(\cdot,\cdot)$. But proving such regularity with the integral $\displaystyle\int_0^t\S_{t-s}\phi_m(x)\d s $ is difficult. So we define 
\begin{align}
\label{3.112} u_m(t,x)=v_m(t,x)-\int_0^t\S_{t-s}\phi_m(x)\d s.
\end{align} This leads to the following regularity result.

\begin{Pro}\label{prop3.6}
Suppose the conditions given in Proposition \ref{prop3.5} are satisfied. Let  $0<\delta<\frac{1}{4}$ and $\gamma<\wi\var_0$. Then there exists a constant $C>0$ and $k(\delta)$ such that 
\begin{align*}
\sup_{t\in[0,T]}t^{(1+\kappa+2\delta)/2}\|\A^{\delta}\D_x u_m(t,\cdot)\|_{0,\gamma+k(\delta)}\leq  C(\gamma,\var,\delta).
\end{align*}
\end{Pro}
\begin{proof}
	Using the exponential estimate given in Proposition \ref{lem3.2} and arguing similar to that of Proposition 4.6, \cite{Da}, one can complete this proof. 
\end{proof}	
\begin{Rem}
Note that the regularity of the  integral $\displaystyle \int_0^t\S_{t-s}\phi_m(s)\d s$ in (\ref{3.112}) cannot be proved as in Proposition \ref{prop3.6}, since 
\begin{align*}
\int_0^t\left\|\A^{\delta}\S_{t-s}\D_x\phi_m(x)\right\|_{0,\gamma+k(\delta)}\d s&\leq \int_0^t\left\|\A^{\delta}\S_{\frac{t-s}{2}}\right\|_{\mathcal{L}(\H)}\left\|\S_{\frac{t-s}{2}}\phi_m(x)\right\|_{1,\gamma+k(\delta)}\d s\no\\&\leq c\int_0^t(t-s)^{-\left(\delta+\frac{(1+\kappa)}{2}\right)}\d s\|\phi_m\|_{0,\gamma+k(\delta)}\no\\&\leq ct^{\frac{1}{2}-\left(\delta+\frac{\kappa}{2}\right)}\|\phi_m\|_{0,\gamma}.
\end{align*}
It is clear from the definition of $\|\cdot\|_{0,\gamma}$ that $\|\phi_m\|_{0,\gamma}$ is not finite, though $t^{\frac{1}{2}-\left(\delta+\frac{\kappa}{2}\right)}$ is bounded for any $0<\delta<\frac{1}{4}$ and $\frac{1}{2}<\kappa<1$.
\end{Rem}
\section{Solvability of Optimal Control Problem}
Now we prove the existence of a mild solution of the HJB equation (\ref{2p16}) using various a-priori estimates of the semigroup, energy estimates derived in previous section, the Arzel\'a-Ascoli theorem and the compactness of $\mathcal{D}(\A^{\delta})$ in $\H$, for  $0<\delta<\frac{1}{4}$. 

In order to prove the boundedness of each term in the HJB equation, we need to bound the integro-differential operator with L\'evy measure in a compact set. It demands  the following assumption on the noise coefficient: 
\begin{Ass}
$	(\A_3)$ For any $x\in\H$ and $0<\theta<1,$  there exists a constant $C>0$ such that   the jump noise coefficient  $\G(\cdot,\cdot)$ satisfies $$\sup_{(t,z)\in[0,T]\times\mZ}\sup_{0<\theta<1}\left\|x+\theta\G(t,z)\right\|\leq C(1+\|x\|).$$
\end{Ass}
Now we are ready to prove the existence of mild solution of the HJB equation (\ref{2p16}) in the smooth function space $\C_{\wi\gamma}^1(\H)$ having exponential weights.

\begin{Thm}\label{thm3.2}
Assume that Assumptions $(\A_1)-(\A_3)$ hold. Then there exists a mild solution $v(t,\cdot)\in\C_{\tilde{\gamma}}^1(\H)$, for some $\tilde{\gamma}\in[0,\wi\var_0]$ where $\wi\var_0=\frac{\var_0}{2}$, of the HJB equation (\ref{2p16}). 
\end{Thm}
\begin{proof}
We first prove the convergence of $u_m(\cdot,\cdot)$ defined in (\ref{3.112}). We extend $u_m(\cdot,\cdot)$ on $[0,T]\times\H$ by setting 
$$u_m(t,x)=u_m(t,\mathrm{P}_mx),\ x\in\H, \ t\in[0,T].$$ By the definition of $u_m(\cdot,\cdot)$, we rewrite the HJB equation as follows: 
\begin{align}\label{eq}
\D_tu_m(t,x)&=\frac{1}{2}\tr(\Q_m \D_x^2u_m)+\left(-\A x+\B_m(x),\D_xu_m\right) + \F(\D_x u_m(t,x))\no \\ &\quad
+\int_{\mathcal{Z}_m}\big(u_m(t,x+\G_m(t,z))-u_m(t,x)-\left(\G_m(t,z),\D_xu_m\right)\big)\mu(\d z),
\end{align}
for any $x\in\mathrm{P}_m\H$ and $t\in[0,T]$. In order to apply the Arzel\'a-Ascoli theorem, we set the compact sets for any $r\in\mathbb{N}$ as follows:  
\begin{equation*}
\left\{
\begin{aligned}
K_r&=\Big\{x\in\mathrm{P}_r\H:\|x\|\leq r\Big\},\ \tau_r=\frac{1}{r},\\
u_m^r&=u_m\Big|_{K_r\times [\tau_r,T]},\ \text{ and }V_m^r=\mathrm{P}_r \D_x u_m\Big|_{K_r\times [\tau_r,T]}=\D_x u_m^r. \end{aligned}
\right.
\end{equation*}
We  now prove the equicontinuity  and boundedness of  $V_m^r(\cdot,\cdot)$. From Proposition (\ref{prop3.5}), we obtain 
\begin{align}\label{3.116}
\sup_{t\in[0,T]}t^{1+\kappa}\|u_m(t,\cdot)\|_{2,\gamma+\var}\leq C(\gamma,\var).
\end{align}
From (\ref{eq}), we have 
\begin{align}\label{3.117}
&\sup_{(t,x)\in K_r\times[\tau_r,T]}|\D_tu_m^r(x,t)|\no\\&\leq \frac{1}{2}\underbrace{\sup_{(x,t)\in K_r\times[\tau_r,T]}|\tr(\Q_m \D_x^2u_m)|}_{:=I_1}+\underbrace{\sup_{(x,t)\in K_r\times[\tau_r,T]}\left|\left(-\A x+\B_m(x),\D_xu_m\right)\right|}_{:=I_2} \no\\&\quad  +\underbrace{\sup_{(x,t)\in K_r\times[\tau_r,T]} |\F(\D_x u_m(t,x))|}_{:=I_3}\no\\&\quad 
+\underbrace{\sup_{(x,t)\in K_r\times[\tau_r,T]}\left|\int_{\mathcal{Z}_m^r}\big(u_m(t,x+\G_m^r(t,z))-u_m(t,x)-\left(\G_m^r(t,z),\D_xu_m\right)\big)\mu(\d z)\right|}_{:=I_4},
\end{align}
where $\G_m^r(\cdot,\cdot)=\mathrm{P}_r\G_m(\cdot,\cdot)$ and $\mathcal{Z}_m^r=\mathrm{P}_r\mathcal{Z}_m$. We have 
\begin{align}\label{3.118}
I_1\leq \sup_{(x,t)\in K_r\times[\tau_r,T]}\tr(\Q_m)| \D_x^2u_m|_{\mathcal{L}^2(\H^2;\mathbb{R})}\leq \tr(\Q)e^{(\gamma+\var)r^2}r^{1+\kappa}\sup_{t\in[0,T]}t^{1+\kappa}\|u_m\|_{2,\gamma+\var}\leq C(\gamma,\var, r).
\end{align}
Since $x=\displaystyle\sum_{i=1}^r(x,e_i)e_i$ and $\lambda_k=\pi^2k^2$, we get \begin{equation}\left\{\begin{aligned}&\|\A x\|^2=\pi^2\sum_{i=1}^r\left\|(x,e_i)\right\|^2i^2\leq \pi^2r^2\|x\|^2\leq \pi^2r^4,\\
&\|\B_m(x)\|^2\leq \|x^2\|_1\leq \|x\|_1^2\leq c\|\A x\|^2\leq c\pi^2 r^4.
\end{aligned}\right.\end{equation}
Using the above estimates, we obtain 
\begin{align}\label{3.119}
I_2&\leq \sup_{(x,t)\in K_r\times[\tau_r,T]}\left(\|\A x\|+\|\B_m(x)\|\right)\|\D_x u_m\|\leq c\pi  e^{(\gamma+\var)r^2}r^{3+\kappa}\sup_{t\in[0,T]}t^{1+\kappa}\|u_m\|_{1,\gamma+\var}\no\\&\leq C(\gamma,\var, r).
\end{align}
From the definition of $\F(\cdot)$, we have
\begin{align}\label{3.120}
I_3\leq \left\{\begin{array}{ll} 
e^{2(\gamma+\var)r^2}r^{2(1+\kappa)}\displaystyle{\sup_{t\in[0,T]}}t^{2(1+\kappa)}\|u_m\|^2_{1,\gamma+\var},
&\|u_m\|_1\leq \rho,\\ \rho e^{(\gamma+\var)r^2}r^{1+\kappa}\displaystyle{\sup_{t\in[0,T]}}t^{1+\kappa}\|u_m\|_{1,\gamma+\var},
&\|u_m\|_1> \rho, \end{array}\right. \leq C(\gamma,\var,r).
\end{align}
For any $x\in K_r,$ by Assumption $(\A_3),$ we have $$\sup_{(t,z)\in[0,T]\times\mZ_m^r}\sup_{0<\theta<1}\left\|x+\theta\G^r_m(t,z)\right\|^2\leq \sup_{(t,z)\in[0,T]\times\mZ}\sup_{0<\theta<1}\|x+\theta\G(t,z)\|^2\leq C (1+r)^2 ,$$
so that  the final term from (\ref{3.117}) can be estimated using  Taylor's theorem   as 
\begin{align}\label{3.121}
I_4&=\frac{1}{2}\sup_{(x,t)\in K_r\times[\tau_r,T]}\left|\int_{\mathcal{Z}_m^r}\D_x^2u_m(t,x+\theta(\G_m^r(t,z))(\G_m^r(t,z),\G_m^r(t,z))\mu(\d z)\right|\nonumber\\&
\leq \frac{1}{2} \sup_{(x,t)\in K_r\times[\tau_r,T]}\left[ \exp\left((\gamma+\var)\sup_{(t,z)\in[0,T]\times\mZ_m^r}\sup_{0<\theta<1}\|x+\theta\G^r_m(t,z)\|^2 \right)t^{-(1+\kappa)}\right]\nonumber\\&
\quad\times\displaystyle{\sup_{t\in[0,T]}}t^{1+\kappa}\|u_m\|_{2,\gamma+\var}\sup_{t\in[0,T]}\int_{{\mZ}}\|\G(t,z)\|^2\mu(\d z)\no
\\&\leq\frac{1}{2} e^{C(\gamma+\var)(1+r)^2}r^{1+\kappa}C(\gamma,\var)\leq C(\gamma,\var,r).
\end{align}
Combining (\ref{3.118})-(\ref{3.121}) and substituting in (\ref{3.117}), we finally obtain 
\begin{align}\label{3.122}
\sup_{(x,t)\in K_r\times[\tau_r,T]}|\D_tu_m^r(x,t)|\leq C(\gamma,\var,r).
\end{align}
Using (\ref{3.116}) and (\ref{3.122}) and arguing as in Theorem 2.2, \cite{Da}, one can prove the equicontinuity of $V_m^r(\cdot,\cdot)$ in $(x,t)\in K_r\times[\tau_r,T]$. Moreover, by Proposition \ref{prop3.6} and Arzel\'a-Ascoli theorem, there exists a continuous function $V^r$ from $K_r\times [\tau_r, T ]$ to $\mathrm{P}_r\H$ such that along a subsequence $$V_{m_k}^r\to V^r\text{ in }\C(K_r\times [\tau_r, T ];\mathrm{P}_r\H).$$ By the compactness of $\D(\A^{\delta})$ in $\H$ and using Proposition \ref{prop3.6}, one can obtain the pointwise strong limit $V(\cdot,\cdot)$ defined on $\displaystyle\left(\bigcup_{r\in\mathbb{N}}\mathrm{P}_r\H\right)\times(0,T]$ of the sequence $V^r(\cdot,\cdot)$. Moreover, $V(\cdot,\cdot)$ is also continuous on the same topology and can be extended to $\H\times(0,T]$. Using this strong convergence of $V^r\to V$ and the estimate given in (\ref{3.116}), and Proposition \ref{prop3.6}, we obtain that $u(\cdot,\cdot)$ is differentiable $\D_xu(x,t)=V(x,t)$. 

In order to prove the convergence of $v_m(\cdot,\cdot)$, we further need the following convergence result whose proof is given at the end of this proof. 

\begin{Lem}\label{lem3.6}
	For any $x\in \H$ and $t\in[0,T]$, we have 
	\begin{align}
	\int_0^t\S_{t-s}^m\phi_m(\mathrm{P}_mx)\d s&\to\int_0^t\S_{t-s}\phi(x)\d s,\label{3.103}\\ \D_x\int_0^t\S_{t-s}^m\phi_m(\mathrm{P}_mx)\d s&\to\D_x\int_0^t\S_{t-s}\phi(x)\d s.\label{3.104}
	\end{align}
\end{Lem}

By the convergence of the sequence $u^m(\cdot,\cdot)$ together with Lemma \ref{lem3.6}, we obtain that $v_{m_k}(t,x)\to v(t,x)$ and 
\begin{align*}
v(t,x)=u(t,x)+\int_0^t\S_{t-s}\phi(s)\d s,
\end{align*}
which is continuous on $[0,T]\times\H$ and $v(0,x)=f(x)$. Moreover, $\D_xv_{m_k}(t,x)\to \D_x v(t,x)$ for any $(t,x)\in(0,T]\times\H$. 

The proof of showing that $v(\cdot,\cdot)$ is indeed the mild solution of (\ref{2p16}) satisfying the integral identity \eqref{2p19a} can be established by arguing similarly as in Theorem 2.2, \cite{Da}.
\end{proof}

\begin{proof}[Proof of Lemma \ref{lem3.6}]
	Let $x\in\H$ and $t\in[0,T]$, then we have 
	\begin{align*}
	&\int_0^t\S_{t-s}^m\phi_m(\mathrm{P}_mx)\d s-\int_0^t\S_{t-s}\phi(x)\d s\no\\
	&=\int_0^t\mE\left[\phi_m(\Y_m(t-s))\right]\d s-\int_0^t\mE\left[\phi(\Y(t-s))\right]\d s\no\\& =\int_0^t\mE\left[\frac{m\|\Y_m(s)\|^2_1}{m+\|\Y_m(s)\|^2}-\|\Y(s)\|^2_1\right]\d s 
	\no\\&\leq \mE\left[\int_0^t\left(\|\Y_m(s)\|_1^2-\|\Y(s)\|_1^2\right)\d s\right]\no\\&\leq \left[\mE\left(\int_0^t\|\Y_m(s)-\Y(s)\|_1^2\d s\right)\right]^{1/2}\left[\mE\left(2\int_0^t\left(\|\Y_m(s)\|_1^2+\|\Y(s)\|_1^2\right)\d s\right)\right]^{1/2}\to 0\text{ as }m\to\infty,
	\end{align*}
	since $\Y_m(\cdot)$ converges to $\Y(\cdot)$ strongly in $\L^2(\Omega;\L^2(0,T;\V))$ and the second integral is bounded by energy estimate in Proposition \ref{prop3.2}. This proves (\ref{3.103}). In order to prove (\ref{3.104}), we first observe that
	\begin{align*}
	&\left(\D_x\int_0^t\S_{t-s}^m\phi_m(\mathrm{P}_mx)\d s-\D_x\int_0^t\S_{t-s}\phi(x)\d s,\mathrm{P}_mh\right)\no\\& = \left(\D_x\int_0^t\mE\left[\phi_m(\Y_m(t-s))\right]\d s-\D_x\int_0^t\mE\left[\phi(\Y(t-s))\right]\d s,\mathrm{P}_mh\right)\no\\& =\int_0^t\left[\mE\left(\D_x\phi_m(\Y_m(t-s)),\eta^{\mathrm{P}_mh}_m(t-s)\right)-\mE\left(\D_x\phi(\Y(t-s)),\eta^{\mathrm{P}_mh}(t-s)\right)\right]\d s\no \\&:=\int_0^t\mE [I(t-s)]\d s.
	\end{align*}
	Next, we write 
	\begin{align*}
	I(s)&=\left(\D_x\phi_m(\Y_m(s))-\D_x\phi_m(\Y(s)),\eta_m^{\mathrm{P}_mh}(s)\right) + \left(\D_x\phi_m(\Y(s))-\D_x\phi(\Y(s)),\eta_m^{\mathrm{P}_mh}(s)\right)\no\\&\quad + \left(\D_x\phi(\Y(s)),\eta_m^{\mathrm{P}_mh}(s)-\eta^{\mathrm{P}_mh}(s)\right):=I_1(s)+I_2(s)+I_3(s).
	\end{align*}
	Note that  
	\begin{align*}
	&\mE\left[\int_0^t|I_1(s)|\d s\right]\no\\&\leq \mE\left[\int_0^t\left|(\D_x\|\Y_m(s)\|_1^2,\eta_m^{\mathrm{P}_mh}(s))\left(1+\frac{\|\Y(s)\|^2}{m}\right)-(\D_x\|\Y(s)\|_1^2,\eta_m^{\mathrm{P}_mh}(s))\left(1+\frac{\|\Y_m(s)\|^2}{m}\right)\right|\d s\right]\no\\&\quad +\frac{1}{m}\mE\left[\int_0^t\|\Y_m(s)\|_1^2|(\D_x\|\Y_m(s)\|^2,\eta_m^{\mathrm{P}_mh}(s))|\d s\right] \no\\&\quad +
\frac{1}{m}\mE\left[\int_0^t\|\Y(s)\|_1^2|(\D_x\|\Y(s)\|^2,\eta_m^{\mathrm{P}_mh}(s))|\d s\right]:=I_4+I_5+I_6. 	\end{align*}
The integral $I_4$ can be further estimated using H\"older's inequality  as 
\begin{align*}
I_4&\leq 2\mE\left[\int_0^t(1+\|\Y(s)\|^2)\|\Y_m(s)-\Y(s)\|_1\|\eta_m^{\mathrm{P}_mh}(s)\|_1\d s\right]\no\\&\quad +\frac{2}{m}\mE\left[\int_0^t\|\Y(s)\|_1\|\eta_m^{\mathrm{P}_mh}(s)\|_1\left(\|\Y(s)\|^2+\|\Y_m(s)\|^2\right)\d s\right]\no\\
&\leq 2\mE\left[\sup_{s\in[0,t]}(1+\|\Y(s)\|^2)\left(\int_0^t\|\Y_m(s)-\Y(s)\|_1^2\d s\right)^{1/2}\left(\int_0^t\|\eta_m^{\mathrm{P}_mh}(s)\|^2_1\d s\right)^{1/2}\right] \no\\&\quad+\frac{2}{m}\mE\left[\sup_{s\in [0,t]}\left(\|\Y(s)\|^2+\|\Y_m(s)\|^2\right)\left(\int_0^t\|\Y(s)\|_1^2\d s\right)^{1/2}\left(\int_0^t\|\eta_m^{\mathrm{P}_mh}(s)\|^2_1\d s\right)^{1/2}\right]\no\\&\leq 2\left\{\mE\left[\sup_{s\in[0,t]}(1+\|\Y(s)\|^2)^4\right]\right\}^{1/4}\left[\mE\left(\int_0^t\|\Y_m(s)-\Y(s)\|_1^2\d s\right)\right]^{1/2}\left\{\mE\left[\left(\int_0^t\|\eta_m^{\mathrm{P}_mh}(s)\|^2_1\d s\right)^2\right]\right\}^{1/4}\no\\&\quad +\frac{2}{m}\left\{\mE\left[\sup_{s\in [0,t]}\left(\|\Y(s)\|^2+\|\Y_m(s)\|^2\right)^4\right]\right\}^{1/4}\left[\mE\left(\int_0^t\|\Y(s)\|_1^2\d s\right)\right]^{1/2}\left\{\mE\left[\left(\int_0^t\|\eta_m^{\mathrm{P}_mh}(s)\|^2_1\d s\right)^2\right]\right\}^{1/4}\no\\&\to 0\ \text{ as }\ m\to\infty.
\end{align*}
In the above calculations, we also used Lemma \ref{lem4.1}, Proposition \ref{prop3.2} and Lemma \ref{lem3.3}. Once again by Proposition \ref{prop3.2} and Lemma \ref{lem4.1}, one can argue that the integrals $I_5$ and $I_6$ also converges to 0 as $m\to\infty$.  Similarly
	\begin{align*}
&	\mE\left[\int_0^t|I_2(s)|\d s\right]\no\\&\leq \mE\left[\int_0^t\left|(\D_x\|\Y(s)\|_1^2,\eta_m^{\mathrm{P}_mh}(s))-(\D_x\|\Y(s)\|_1^2,\eta_m^{\mathrm{P}_mh}(s))\left(1+\frac{\|\Y(s)\|^2}{m}\right)\right|\d s\right]\no\\&\quad +\frac{1}{m}\mE\left[\int_0^t\|\Y(s)\|_1^2|(\D_x\|\Y(s)\|_1^2,\eta_m^{\mathrm{P}_mh}(s))|\d s\right]\no\\&\leq \frac{2}{m}\mE\left[\int_0^t\|\Y(s)\|_1\|\eta_m^{\mathrm{P}_mh}(s)\|_1\|\Y(s)\|^2\d s\right] + \frac{2}{m}\mE\left[\int_0^t\|\Y(s)\|_1^2\|\eta_m^{\mathrm{P}_mh}(s)\|\|\Y(s)\|\d s\right]\no\\&\leq \frac{2}{m}\left[\mE\left(\sup_{s\in [0,t]}\|\Y(s)\|^4\right)\right]^{1/4}\left[\mE\left(\int_0^t\|\Y(s)\|_1^2\d s\right)\right]^{1/2}\left\{\mE\left[\left(\int_0^t\|\eta_m^{\mathrm{P}_mh}(s)\|_1^2\d s\right)^2\right]\right\}^{1/4}\no\\&\quad+\frac{2}{m}\left[\mE\left(\sup_{s\in [0,t]}\|\Y(s)\|^4\right)\right]^{1/4}\left[\mE\left(\sup_{s\in [0,t]}\|\eta_m^{\mathrm{P}_mh}(s)\|^4\right)\right]^{1/4}\left\{\mE\left[\left(\int_0^t\|\Y(s)\|_1^2\d s\right)^2\right]\right\}^{1/2}\no\\&\to 0\ \text{ as } \ m\to\infty.
	\end{align*}
	An application of the convergence (\ref{3.77}) and the energy estimates lead to the convergence of $\displaystyle\mE\left[\int_0^t|I_3(s)|\d s\right]$ to $0$ as $m\to\infty$. This completes the proof.
\end{proof}

Now using the solution $v_m\in \C_{\wi\gamma}^1(\H)$, we justify the feedback control formula $\wi\U_m(t)=\mathcal{G}( \D_xv_m(T-t,\wi \X_m(t)))$ and prove the convergence of $\wi\X_m(\cdot)$  to $\wi\X(\cdot)$. Moreover, we establish an identity satisfied by the cost functional $\mathcal{J}(0,T;x,\U)$.

\begin{Thm}\label{thm3.3}
Suppose the conditions given in Theorem \ref{thm3.2} are satisfied. Then for any control $\U\in\mathscr{U}_{\rho}^{0,T}$, the following identity holds:
\begin{align}\label{3.122a}
\mathcal{J}(0,T;x,\U)=v(T,x)+\frac{1}{2}\mE\left[\int_0^T\|\U(t)+\D_x v(T-t,\X(t))\|^2-\chi(\|\D_x v(T-t,\X(t))\|-\rho)\d t\right],
\end{align}
where $\mathcal{J}(0,T;x,\U)$ is defined in (\ref{2p13b}), the function $\chi$ satisfies  $\chi(a)=0$ for $a\leq 0$ and $\chi(a)=a^2$ for $a\geq 0$ and $\X(\cdot)$ is the solution of (\ref{2p12}). Moreover, the closed loop equation (\ref{2p18}) has an  optimal pair  $(\wi\X,\wi\U)$ with $\wi \U(t)= {\mathcal G}( \D_xv(T-t,\wi \X(t)))$.
\end{Thm}
In order to prove Theorem \ref{thm3.3}, we need the following results on the approximated cost functional. 

\begin{Lem}
For $m\in\mathbb{N}$, let $\X_m(\cdot)$ be the solution of (\ref{2F1}) and $\U_m =\mathrm{P}_m\U\in \L^2(\Omega;\L^2(0,T;\mathrm{P}_m\H))$. For any $x\in\mathrm{P}_m\H$, the following holds for the approximated cost functional  
\begin{align*}
&\mathcal{J}_m(0,T;x,\U_m)\no\\&=v_m(T,x)+\frac{1}{2}\mE\left[\int_0^T\Big(\|\U_m(t)+\D_xv_m(T-t,\X_m(t))\|^2-\chi(\|\D_x v_m(T-t,\X_m(t))\|-\rho)\Big)\d t\right], 
\end{align*}
where $\chi(\cdot)$ is the function defined in Theorem \ref{thm3.3}.
\end{Lem}
\begin{proof}
Let us apply the finite dimensional It\^o formula to the process $v_m(t-s,\X_m(s))$ to obtain 
\begin{align}\label{3.124}
 \d v_m(t-s,\X_m(s))&=\Big[-\D_sv_m(t-s,\X_m(s))+\mathscr{L}_xv_m(t-s,\X_m(s))\no\\&\qquad +(\D_xv_m(t-s,\X_m(s)),\U_m(s))\Big]\d s+\d\M^m_{s}\no\\&= \Big[-\phi_m(\X_m(s))-\F_m(\D_xv_m(t-s,\X_m(s)))\no\\ &\qquad+(\D_xv_m(t-s,\X_m(s)),\U_m(s))\Big]\d s +\d\M^m_{s},
\end{align}
where 
\begin{align*}
\d\M^m_{s}&=(\D_xv_m(t-s,\X_m(s)),\Q_m^{1/2}\d\W_m(s)) \no\\&\quad+\int_{{\mZ}_m}\left[v_m(t-s,\X_m(s)+\G_m(s,z))-v_m(t-s,\X_m(s))\right]\wi{\mathrm{N}}(\d s,\d z).
\end{align*}
Setting $t=T$, integrating from $0$ to $T$, and taking expectation, we get 
\begin{align}\label{3.125}
v_m(0,\X_m(T))&=v_m(T,x)+\mE\left[\int_0^T(\D_xv_m(T-s,\X_m(s)),\U_m(s))\d s\right]\no\\&\quad -\mE\left[\int_0^T\Big[\phi_m(\X_m(s))+\F_m(\D_xv_m(T-s,\X_m(s)))\Big]\d s\right].
\end{align} 
Using the definition of approximated cost functional $\mathcal{J}_m$ and arranging the terms in (\ref{3.125}), we arrive at 
\begin{align}\label{3.126}
\mathcal{J}_m(0,T;x,\U_m)&= \mE\left[\int_0^T\left(\phi_m(\X_m(t))+\frac{1}{2}\|\U_m(t)\|^2\right)\d t+f_m(\X_m(T))\right]\no\\&=v_m(T,x)+\frac{1}{2}\mE\left[\int_0^T\|\U_m(t)+\D_xv_m(T-t,\X_m(t))\|^2\d t\right]\no\\&\quad -\mE\left[\int_0^T\Big(\F_m(\D_xv_m(T-t,\X_m(t)))+\frac{1}{2}\|\D_xv_m(T-t,\X_m(t))\|^2\Big)\d t\right].
\end{align}
Using the definition of $\F_m(\cdot)$, the final integral in (\ref{3.126}) can be written as 
\begin{align*}
\frac{1}{2}\mE\left[\int_0^T\chi(\|\D_x v_m(T-t,\X(t))\|-\rho)\d t\right].
\end{align*}
This completes the proof.
\end{proof}

\begin{Lem}
The approximated solution  $\X_m(\cdot)$ of (\ref{2F1}) converges almost surely to the solution $\X(\cdot)$ of (\ref{2p12}) in $\L^2(0,T;\V)\cap\mathscr{D}([0,T];\H)$. Moreover, the approximated cost functional $\mathcal{J}_m(0,T;x,\U_m)$ converges  to $\mathcal{J}(0,T;x,\U)$.
\end{Lem}
\begin{proof}
It is clear from (\ref{c1}) and the rest of the arguments in Proposition \ref{prop3.3} that in order  to prove the almost sure convergence of $\X_m(\cdot)$ in the given topology, it is sufficient to prove 
\begin{align}\label{3.127}
\mE\left[\int_0^T(\U_m(t),\X_m(t))\d t\right]\to\mE\left[\int_0^T(\U(t),\X(t))\d t\right], \text{ as } m\to\infty.
\end{align}
Since $\X_m\rightharpoonup\X$ in $\L^2(\Omega;\L^2(0,T;\V))$ and $\U_m\to\U$ in $\L^2(\Omega;\L^2(0,T;\H))$, the convergence in (\ref{3.127}) follows. 

Since $\|\U_m(t)\|\leq \rho$, $\mP-$a.s., by  Proposition \ref{lem3.2} with   $\displaystyle 0<\var < \widehat\var_0(\var_0)$ leads to
\begin{align}\label{3.128}
\mE\left[\exp\left(\var \|\X_m(t)\|^2 + \var \int _0^t \|\X_m(s)\|^2_1 \d s\right)\right] 
\leq C(\var,\Q,T,\rho)  \exp(\var \|x\|^2).
\end{align}
Along with the estimate (\ref{3.128}), using similar arguments as in  the proof of Theorem 2.3, \cite{Da}, one can prove that $\mathcal{J}_m(0,T;x,\U_m)$ converges to $\mathcal{J}(0,T;x,\U)$ satisfying (\ref{3.122a}).
\end{proof}

\begin{proof}[Proof of Theorem \ref{thm3.3}] 
	It can be easily shown that the following approximated closed loop equation  has a unique solution $\wi\X_m(\cdot)$:
	\begin{equation} 
	\left\{\begin{aligned}
	\d\wi \X_m(t)&=[-\A\wi \X_m(t)+\B(\wi \X_m(t))+{\mathcal G}(\D_xv(T-t,\wi \X_m(t)))]\d t+\Q_m^{1/2}\d\W_m(t)\\
	&\quad+ \int_{\mZ_m}\G_m(t,z) \wi{\mathrm{N}}(\d t,\d z), \;t\in (0,T), \\
	\wi  \X_m(0)&=\mathrm{P}_mx.
	\end{aligned}
	\right.
	\end{equation}
From the definition of $\mathcal{G}$ in (\ref{2p17}), we know that $$\|\mathcal{G}(p)\|\leq \rho,\ \text{ for all }\ p\in\H.$$ Arguing similarly as in Proposition \ref{prop3.3}, one can obtain that $\wi\X_m(\cdot)$ is almost surely bounded in $\mathscr{D}([0,T];\H)\cap\L^2(0,T;\V).$  Then, there exists a subsequence converging to the solution $\wi\X(\cdot)$ of the closed loop equation (\ref{2p18}). If we show that the closed loop equation (\ref{2p18}) has at most one solution, the above convergence along the subsequence imply that the entire sequence converges.

Now, we prove the uniqueness of the closed loop equation (\ref{2p18}). Let $\X_1(\cdot)$ and $\X_2(\cdot)$ be two solution of (\ref{2p18}) and $\wi\X(\cdot)=\wi\X_1(\cdot)-\wi\X_2(\cdot)$. Then, we have 
\begin{equation} \label{3.131}
\frac{\d}{\d t}\wi \X(t)=-\A\wi \X(t)+\B(\wi \X_1(t))-\B(\wi \X_2(t))+{\mathcal G}(\D_xv(T-t,\wi \X_1(t)))-{\mathcal G}(\D_xv(T-t,\wi \X_2(t))).
\end{equation}
Taking inner product with $\wi\X(\cdot)$ in (\ref{3.131}), we get
\begin{align}\label{3.132}
\frac{1}{2}\frac{\d}{\d t}\|\wi\X(t)\|^2+\|\wi\X(t)\|_1^2&=\left(\B(\wi \X_1(t))-\B(\wi \X_2(t)),\wi\X(t)\right)\no\\&\quad+\left({\mathcal G}(\D_xv(T-t,\wi \X_1(t)))-{\mathcal G}(\D_xv(T-t,\wi \X_2(t))),\wi\X(t)\right)\no\\&:=I_1+I_2.
\end{align}
Using the Gagliardo-Nirenberg inequality, the nonlinear term $I_1$ can be estimates as 
\begin{align*}
\left|\left(\B(\wi \X_1)-\B(\wi \X_2),\wi\X\right)\right|&\leq C\|\wi \X\|^{3/2}\|\wi \X\|_1^{1/2}\left(\|\wi \X_1\|_1+\|\wi \X_2\|_1\right)\no\\&\leq \|\wi \X\|_1^2+C\|\wi \X\|^2\left(\|\wi \X_1\|_1^{4/3}+\|\wi \X_2\|_1^{4/3}\right).
\end{align*}
From the definition of $\mathcal{G}(\cdot)$, it is clear that
\begin{align}\label{3.134}
|I_2|&\leq \|{\mathcal G}(\D_xv(T-t,\wi \X_1))-{\mathcal G}(\D_xv(T-t,\wi \X_2))\| \|\wi\X\|\no\\&\leq \|\D_xv(T-t,\wi \X_1)-\D_xv(T-t,\wi \X_2)\|\|\wi\X\|\no\\&=\left\|\int_0^1\frac{\d}{\d r}\D_x v(T-t,r\wi\X_1+(1-r)\wi\X_2)\d r\right\|\|\wi\X\|\leq   |\D_x^2v(T-t,y)|_{\mathcal{L}^2(\H^2;\mR)}\|\wi\X\|^2\no\\&\leq C e^{(\gamma+\var)M}(T-t)^{-(1+\kappa)}\sup_{t\in[0,T]}(T-t)^{(1+\kappa)}\|v(T-t)\|_{2,\gamma+\var}\|\wi\X\|^2\no\\&\leq C(\gamma,\var)e^{(\gamma+\var)M}(T-t)^{-(1+\kappa)}\|\wi\X\|^2,
\end{align}
 for $\|\D_xv(T-t,\wi \X_1)\|\leq \rho$ and $\|\D_xv(T-t,\wi \X_2)\|\leq \rho$, where $\displaystyle M=\max\left\{\sup_{t\in[0,T]}\|\wi\X_1\|^2,\sup_{t\in[0,T]}\|\wi\X_2\|^2\right\}$. For the last inequality in (\ref{3.134}), we also used Proposition \ref{prop3.5}. For $\|\D_xv(T-t,\wi \X_1)\|>  \rho$ and $\|\D_xv(T-t,\wi \X_2)\|> \rho$, we have 
 \begin{align}\label{3.134a}
 |I_2|&\leq \left\|-\rho\frac{\D_xv(T-t,\wi \X_1)}{\|\D_xv(T-t,\wi \X_1)\|}+\rho\frac{\D_xv(T-t,\wi \X_2)}{\|\D_xv(T-t,\wi \X_2)\|}\right\|\|\wi\X\|\no\\&\leq \frac{\rho}{\|\D_xv(T-t,\wi \X_1)\|\|\D_xv(T-t,\wi \X_2)\|}\left\|\left(\|\D_xv(T-t,\wi \X_1)\|-\|\D_xv(T-t,\wi \X_2)\|\right)\D_xv(T-t,\wi \X_1)\right.\no\\&\left.\quad +\|\D_xv(T-t,\wi \X_1)\|\left(\D_xv(T-t,\wi \X_2)-\D_xv(T-t,\wi \X_1)\right)\right\|\|\wi\X\|\no\\&\leq \frac{2\rho}{\|\D_xv(T-t,\wi \X_2)\|}\|\D_xv(T-t,\wi \X_1)-\D_xv(T-t,\wi \X_2)\|\|\wi\X\|\no\\&\leq C(\gamma,\var)e^{(\gamma+\var)M}(T-t)^{-(1+\kappa)}\|\wi\X\|^2,
 \end{align}
 using \eqref{3.134}. Now for $\|\D_xv(T-t,\wi \X_1)\|>  \rho$ and $\|\D_xv(T-t,\wi \X_2)\|\leq  \rho$, we obtain 
 \begin{align}\label{3.134b}
 |I_2|&\leq \left\|-\rho\frac{\D_xv(T-t,\wi \X_1)}{\|\D_xv(T-t,\wi \X_1)\|}+\D_xv(T-t,\wi \X_2)\right\|\|\wi\X\|\no\\&\leq \frac{1}{\|\D_xv(T-t,\wi \X_1)\|}\left\| \|\D_xv(T-t,\wi \X_1)\|\left(\D_xv(T-t,\wi \X_2)-\D_xv(T-t,\wi \X_1)\right)\right.\no\\&\left.\quad +\D_xv(T-t,\wi \X_1)\left(\|\D_xv(T-t,\wi \X_1)\|-\rho\right)\right\| \|\wi\X\|\no\\&\leq \left(\|\D_xv(T-t,\wi \X_1)-\D_xv(T-t,\wi \X_2)\|+\|\D_xv(T-t,\wi \X_1)\|-\rho\right)\|\wi\X\|\no\\&\leq C(\gamma,\var)e^{(\gamma+\var)M}(T-t)^{-(1+\kappa)}\|\wi\X\|^2,
 \end{align}
 using \eqref{3.134}. An estimate similar to \eqref{3.134b} can be obtained for $\|\D_xv(T-t,\wi \X_1)\|\leq  \rho$ and $\|\D_xv(T-t,\wi \X_2)\|>  \rho$. Now the equation (\ref{3.132}) can be estimated by the Gr\"onwall inequality as follows: 
\begin{align*}
\|\wi\X(t)\|^2\leq \|\wi\X(0)\|^2\exp\left(\int_0^tC\left(\|\wi \X_1(s)\|_1^{4/3}+\|\wi \X_2(s)\|_1^{4/3}\right)\d t+C(\gamma,\var,M)(T-t)^{-\kappa}\right),
\end{align*}
for any $t\in[0,T)$. By an estimation similar to \eqref{Y1} and Proposition \ref{lem3.2}, the exponential integral is bounded.  Thus, if $\wi\X_1(0)=\wi\X_2(0)=x$, then $\wi\X_1(t)=\wi\X_2(t)$ for all $t\in[0,T)$, $\mP-$a.s. If there is a jump at $t=T$, we can extend the analysis for some time $\wi T>T$ and can conclude the uniqueness for all $t\in[0,\wi T]$. Otherwise the uniqueness follows immediately. 

Moreover, $\wi\U(\cdot)$ is an optimal control. Indeed, from (\ref{3.122a}), it is clear that $$v(T,x)\leq \mathcal{J}(0,T;x,\U),\text{ for all }\U\in\L^2(\Omega;\L^2(0,T;\H)).$$ Since $\wi\X(\cdot)$ is the solution corresponding to the control $\wi\U(t)={\mathcal G}( \D_xv(T-t,\wi \X(t)))$, from (\ref{3.122a}) we also have $v(T,x)=\mathcal{J}(0,T;x,\wi\U)$. From the above inequality it is clear that $\wi\U(\cdot)$ is an optimal control and $(\wi\X,\wi\U)$ is an optimal pair.
\end{proof}

	\medskip\noindent
{\bf Acknowledgements:} M. T. Mohan would  like to thank the Department of Science and Technology (DST), India for Innovation in Science Pursuit for Inspired Research (INSPIRE) Faculty Award (IFA17-MA110), and Indian Statistical Institute-Bangalore Centre (ISI-Bangalore), Indian Institute of Technology Roorkee (IIT-Roorkee) and  Indian Institute of Space Science and Technology (IIST), Trivandrum for providing stimulating scientific environment and resources. S. S. Sritharan's work has been funded by U. S. Army Research Office, Probability and Statistics program.  Authors would also like to thank Prof. A. Debussche for useful suggestions in proving Theorem  \ref{thm3.2}.  The authors sincerely would like to thank the reviewers for their  valuable comments and suggestions, which led to the improvement  of this paper.

\begin{appendix}
	\renewcommand{\thesection}{\Alph{section}}
	\numberwithin{equation}{section}
\section{Dynamic Programming Principle} \label{DPP}
\setcounter{section}{1}
\setcounter{equation}{0}
 \renewcommand{\thesection}{\Alph{section}} 

We briefly give  the key steps in deriving the stochastic Dynamic Programming Principle and  HJB equation.  Let us define 
\begin{align*}
\mathfrak{L}(\X,\U)=\|\D_\xi \X\|^2+\frac{1}{2}\|\U\|^2, \ \   \Psi(\X(T))=\|\X(T)\|^2, \ \
\F(\X,\U)=-\A\X+\B(\X)+\U.
\end{align*} 
Let $T>0$ be given. Then for  any $t\in[0,T),$ consider the problem of minimizing the cost functional 
\begin{align} \label{a1}
{\mathcal J}(t,T;x,\U)= \mE_t\left[\int_t^T\mathfrak{L}(\X(s),\U(s))ds+ \Psi(\X(T)) \right],
\end{align}
over all controls $\U\in \mU^{t,T}_{\rho}$  and  $\X(\cdot)$ satisfying   the following state equation  
\begin{equation} \label{a2}
\left\{
\begin{aligned}
\d  \X(s)&=\F(\X(s),\U(s))\d s+\Q^{1/2}\d\W(s) +\int_{\mZ}\G(s,z) \wi{\mathrm{N}}(\d s,\d z),  \ s\in [t,T], \\
\X(t)&=x, \ x\in \H.
\end{aligned} 
\right.
\end{equation}
 Here $\mE_t[\X(s)]=\mE[\X(s) | \X(t)=x]$, for $s\geq t$  and the admissible control set is defined as
\begin{align*}
\mU^{t,T}_{\rho} =\Big\{\U\in \L^2(\Omega, \L^2(t,T;\H)):\|\U\|\leq\rho, \ \mbox{and} \  \U \ \mbox{is adapted to} \ \mF_{t,s}\Big\},
\end{align*}
where $\mF_{t,s}$ is the $\sigma$-algebra  generated by the paths of $\W$ and random measures $\mathrm{N}$ upto time $s,$ i.e.,  $\sigma\{\W(r); t\leq r\leq s\}$ and $\sigma\{\mathrm{N}(S); S\in \mB([t,s]\times\mZ)\}.$      The value function of the above control problem is defined as  
\begin{equation*}
\left\{
\begin{aligned}
\mathscr{V}(t,x)&=\inf_{\U\in\mU^{t,T}_{\rho}} {\mathcal J}(t,T;x,\U),   \text{ for all }t\in[0,T), \ x\in \H, \\ 
\mathscr{V}(T,x)&=\|x\|^2.
\end{aligned} 
\right.
\end{equation*}
Following the dynamic programming strategy,  any admissible control  $\U\in\mU^{t,T}_{\rho}$ is the combination of  controls in $ \mU^{t,\tau}_{\rho}$ and $ \mU^{\tau,T}_{\rho}$ for any $0\leq t< \tau< T.$ More precisely,  suppose the processes $\U_1(s)$ and $\U_2(s)$ be the restriction of the control process $\U(s)$ to the time intervals $[t,\tau] $ and $[\tau,T] $  respectively, i.e.,
\begin{equation*}
\U(s)= (\U_1\oplus\U_2)(s)
=\left\{
\begin{array}{ll}
\U_t(s),& s\in[t,\tau], \\
\U_\tau(s),& s\in[\tau, T]. \\
\end{array}\right.
\end{equation*}
Accordingly, the  admissible control set is written as $\mU^{t,T}_{\rho}=\mU^{t,\tau}_{\rho}\oplus \mU^{\tau,T}_{\rho}.$ Note also that the controls $\U_1(s)$ and $ \U_2(s)$ are adapted to $\mF_{t,s}$ and $\mF_{\tau,s}$ respectively. 

The system state $\X(\cdot)$ is determined by \eqref{a2} with $\U(s)=(\U_1\oplus\U_2)(s)\in \mU^{t,T}_{\rho}.$  We decompose the system state as $\X(s)=(\X_1\oplus  \X_2)(s),$ where $\X_1$ and $\X_2$ satisfy 
\begin{equation}\label{a4}
\left\{
\begin{aligned}
\d  \X_1(s)&=\F(\X_1(s),\U_1(s))\d s+\Q^{1/2}\d\W(s) +\int_{\mZ}\G(s,z) \wi{\mathrm{N}}(\d s,\d z),  \ s\in [t,\tau], \\
\X_1(t)&=\X(t)=x\in \H, 
\end{aligned}
\right.
\end{equation}
and 
\begin{equation} \label{a5}
\left\{
\begin{aligned}
\d  \X_2(s)&=\F(\X_2(s),\U_2(s))\d s+\Q^{1/2}\d\W(s) +\int_{\mZ}\G(s,z) \wi{\mathrm{N}}(\d s,\d z),  \ s\in [\tau,T],\\
\X_2(\tau)&=\X_1(\tau)=\X(\tau).
\end{aligned}
\right.
\end{equation} 
By the tower property of conditional expectation 
$$\mE_t\Big[\mE_t\big(\mathfrak{L}(\X_2(s),\U_2(s))\ |  \mF_{t,\tau}\big)\Big] =\mE_t\big[\mathfrak{L}(\X_2(s),\U_2(s))\big], \tau\leq s\leq T, $$   we write
\begin{align*} 
\mathscr{V}(t,x)&= \inf_{\U\in\mU^{t,T}_{\rho} }\mE_t\left\{\int_t^\tau\mathfrak{L}(\X(s),\U(s))\d s + \int_\tau^T\mathfrak{L}(\X(s),\U(s))\d s+\Psi(\X(T))\right\}  \\ 
&= \inf_{\substack{\U_1\in\mU^{t,\tau}_{\rho}, \U_2\in\mU^{\tau,T}_{\rho} \\  \X_2(\tau)=\X_1(\tau) }} \mE_t\left\{ \int_t^\tau\mathfrak{L}(\X_1(s),\U_1(s))\d s \right.\\
&\left. \quad+ \mE_t\left[ \int_\tau^T\mathfrak{L}(\X_2(s),\U_2(s))\d s+\Psi(\X_2(T))\big|\mF_{t,\tau}\right]\right \}. 
\end{align*}
Using the Markovian property of the process $\X(\cdot)$ one can get 
for $\tau\leq s\leq T,$   $$\mE_t\Big[\mathfrak{L}(\X_2(s),\U_2(s))\ |  \mF_{t,\tau}\Big] =\mE_{\tau}\Big[\mathfrak{L}(\X_2(s),\U_2(s))\ | \ \X_2(\tau)=\X(\tau) \Big]$$ and the same reasoning is true for $\Psi(\X_2(T))$ as well. It leads to  
\begin{align} \label{a21}
{\mathcal J}(\tau,T;\X_2(\tau),\U_2)= \mE_t\left\{\int_\tau^T\mathfrak{L}(\X_2(s),\U_2(s))\d s+\Psi(\X_2(T)) \ \Big | \mF_{t,\tau}\right\}.
\end{align} Besides, for more details on the case of continuous diffusion, one can  refer to \cite{Y} (and also \cite{G}) and that can be modified to this case. Hence
\begin{align}\label{a3}
\mathscr{V}(t,x)&=  \inf_{\U_1\in\mU^{t,\tau}_{\rho}}  \mE_t\left\{\int_t^\tau\mathfrak{L}(\X_1(s),\U_1(s))\d s\right\}\no  \\
&\quad+ \mE_t\left\{\inf_{\substack{\U_2\in\mU^{\tau,T}_{\rho} \\ \X_2(\tau)=\X_1(\tau) }}\mE_t\left[\int_\tau^T\mathfrak{L}(\X_2(s),\U_2(s))\d s+\Psi(\X_2(T))\Big| \mF_{t,\tau}\right]\right\} \no \\ 
&= \inf_{\U_1\in\mU^{t,\tau}_{\rho}}  \mE_t\left\{\int_t^\tau\mathfrak{L}(\X_1(s),\U_1(s))\d s  +\mathscr{V}(\tau,\X_1(\tau)) \right\}. 
\end{align}
Thus, we have proved the following  DPP {(or \it Bellman's principle of optimality)}
\begin{align} \label{a6}
\mathscr{V}(t,x)=\inf_{\U\in\mU^{t,\tau}_{\rho}}  \mE\left\{\int_t^\tau\mathfrak{L}(\X(s),\U(s))\d s  +\mathscr{V}(\tau,\X(\tau)) \big| \X(t)=x\right\}.
\end{align}

\subsection{Dynamic Programming Equation} 

Rewrite \eqref{a6} as follows
\begin{align}\label{a7}
\inf_{\U\in\mU^{t,\tau}_{\rho}}  \mE\left\{\int_t^\tau\mathfrak{L}(\X(s),\U(s))ds  +\mathscr{V}(\tau,\X(\tau)) - \mathscr{V}(t,x) \big| \X(t)=x\right\}=0.
\end{align}
Now the HJB equation can be formally obtained by applying the  It\^o formula for $\mathscr{V}(\tau,\X(\tau)) - \mathscr{V}(t,x).$ 

\begin{Rem}[It\^o Formula]
For any $t\geq 0$ and   $\Phi\in \C_b^{1,2}([0,T];\mathcal{D}(\mathscr{A}^{\U})),$ the following relation holds (in fact,  $\Phi$ could be a function such that whose Gate\^aux derivatives $\D_x\Phi$ and $\D^2_x\Phi$ are H\"older continuous with  exponent $\delta=1$, see, \cite{Me}): 
\begin{align} \label{a8}
\Phi(\tau,\X(\tau))-\Phi(t,x)=\int_t^\tau\big[\D_s\Phi(s,\X(s))+\mA^{\U}\Phi(s,\X(s))\big]\d s + \M_\tau,  \ \mP{\text -a.s.},
\end{align}
where $\mA^{\U}\Phi$ is the second order partial integrodifferential operator
\begin{align*}
\mA^{\U}\Phi(s,\X(s))& = \la \D_x\Phi(s,\X(s)), \F(\X(s),\U(s))\ra +\frac{1}{2} \tr(\Q \D^2_x\Phi(s,\X(s)) \\
&\quad+\iZ\big[\Phi(s,\X(s)+\G(s,z))-\Phi(s,\X(s))-\la \G(s,z),\D_x\Phi(s,\X(s))\ra\big]\mu(\d z)
\end{align*}
and $\M_{\tau}$ is the martingale given by
\begin{align*}
\M_\tau=&\int_t^\tau\la  \D_x\Phi(s,\X(s)), \sqrt{\Q}\d\W(s)\ra \\
&+ \int_t^\tau\iZ \big[\Phi(s,\X(s-)+\G(s,z))-\Phi(s,\X(s-))\big]\wi{\mathrm{N}}(\d s,\d z).
\end{align*}
\end{Rem}
Taking  $\Phi(\tau,\X(\tau))=\mathscr{V}(\tau,\X(\tau))$ (of course by assuming required smoothness and boundedness on $\mathscr{V}$)  and plug this back into \eqref{a7}, and use the fact that martingale $\M_{\tau}$ has a zero mean to obtain the following
\begin{equation}\label{a9}
\inf_{\U\in\mU^{t,\tau}_{\rho}}  \mE\left\{\int_t^\tau\big[\mathfrak{L}(\X(s),\U(s))+\D_s\mathscr{V}(s,\X(s))+\mA^{\U}\mathscr{V}(s,\X(s))\big]\d s \big| \X(t)=x\right\}=0.
\end{equation}
Finally, take $\tau=t+h, h>0$ and divide by $h.$  Passing $h\to 0$ and using the conditional expectation, we formally obtain the HJB equation 
\begin{equation}\label{a10}
\left\{
\begin{aligned}
&\D_t\mathscr{V}(t,x)+\inf_{\|\U\|\leq \rho} \big\{\mA^{\U}\mathscr{V}(t,x) +\mathfrak{L}(x,\U(t)) \big\}=0, \ \ t\in[0,T), \\
&\mathscr{V}(T,x)=\|x\|^2,  \ x\in \H.
\end{aligned}
\right.
\end{equation}
Moreover, setting $v(t,x)=\mathscr{V}(T-t,x),$ one can obtain the following initial value problem
\begin{equation}\label{a11}
\left\{
\begin{aligned} 
\D_tv(t,x)=&  \mathscr{H}(x,t,\D_xv,\D^2_xv), \ \  \ t\in (0,T), \\
  v(0,x)=& \|x\|^2, \ x\in \H, \no
\end{aligned} 
\right.\end{equation}
where $\mathscr{H}$ is given by
\begin{align*}
\mathscr{H}(x,t,\D_xv,\D^2_xv) =&\frac{1}{2}\tr(\Q \D_x^2v)+\la - \A x+\B(x),\D_xv\ra  \\
&+\iZ\big[v(t,x+\G(t,z))-v(t,x)-\la \G(t,z),\D_xv\ra\big]\mu(\d z) \\
&+\|\D_\xi x\|^2+\inf_{\|\U\|\leq \rho}\left\{\la \U(t),\D_xv\ra+\frac{1}{2}\|\U(t)\|^2\right\}.
\end{align*}  

\end{appendix}



{\bf Data Sharing Statement}: Data sharing not applicable to this article as no datasets were generated or analysed during the current study.


\end{document}